\documentclass[12pt,leqno]{amsart}

\usepackage{enumerate} 
\usepackage[utf8]{inputenc} 
\usepackage[T1]{fontenc} 
\usepackage{amssymb,amsthm}
\usepackage{mathrsfs}
\usepackage{textcomp}
\usepackage{bm}
\usepackage{graphicx}
\usepackage{xcolor}

\newtheorem{thm}{Theorem}[section]
\newtheorem{prop}[thm]{Proposition}
\newtheorem{lem}[thm]{Lemma}
\newtheorem{cor}[thm]{Corollary}

\theoremstyle{definition}

\theoremstyle{remark}

\numberwithin{equation}{section}

\usepackage[font=footnotesize,tableposition=top,figureposition=bottom,skip=0pt]{caption}
\usepackage{hyperref}
\hypersetup{%
	linkcolor=blue,
	anchorcolor=black,
	citecolor=red,
	filecolor=magenta,
	menucolor=red,
	urlcolor=magenta,
	colorlinks=true,
	bookmarksnumbered=true,
	bookmarksopen=true,
	bookmarksopenlevel=0,
	breaklinks=true
}

\begin{document}

\title[Trigonometric chaos and Random distributions]
{Trigonometric multiplicative chaos and  Application to random distributions  \footnote{This paper has been accepted for publication in SCIENCE CHINA Mathematics.} }

\author{Ai-hua Fan}
\address{School of mathematics and Statistics, Central China Normal University,  430079 Wuhan, China \& LAMFA, UMR 7352 CNRS, University of Picardie, 33 rue Saint Leu,80039 Amiens, France}
\email{ai-hua.fan@u-picardie.fr}
\author{Yves Meyer }
\address{CMLA, ENS-Cachan, CNRS, Universit\'e Paris-Saclay
France}
\email{yves.meyer@ens-paris-saclay.fr}
\thanks{}
\maketitle






\begin{abstract}
The random trigonometric series  $\sum_{n=1}^\infty \rho_n \cos (nt +\omega_n)$ on the circle  $\mathbb{T}$ are studied under the conditions
$\sum |\rho_n|^2=\infty$ and $\rho_n\to 0$, where $\{\omega_n\}$ are iid and uniformly distributed on $\mathbb{T}$. They are almost surely not  Fourier-Stieljes series but define pseudo-functions.
This leads us to develop the theory of trigonometric multiplicative chaos,
which produces a class of random measures.  The behaviors of the partial sums of the above series  are proved to be multifractal.
Our theory holds on the torus $\mathbb{T}^d$ of dimension $d\ge 1$.

\end{abstract}

{\small	
  \textbf{Keywords--}   Multiplicative chaos, Random Fourier series, Hausdorff dimension, Riesz potential
}

{\small	
  \textbf{MSC(2020)--}   60G57, 60G46, 42A20,28A78
}

\tableofcontents

\section{Introduction}

Let us consider a random trigonometric series 
\begin{equation}\label{eq:RTS}
 \sum_{n=1}^\infty X_n \cos  (nt +\Phi_n)
\end{equation} 
where the random complex variables $A_n:= X_n e^{i\Phi_n}$ ($X_n, \Phi_n$ being real) are supposed to be independent and symmetric. It is the real part of the 
random Taylor trigonometric series 
\begin{equation}\label{eq:RTS2}
 \sum_{n=1}^\infty A_n e^{int}.
\end{equation} 
These series are well studied  in \cite{PZ1930-1932} and \cite{Kahane1960} (see also \cite{Kahane1985}) under the condition $\mathbb{E}X_n^4 \le C (\mathbb{E}X_n^2)^2$. 
If $\sum \mathbb{E}(X_n^2) =+\infty $, the series (\ref{eq:RTS}) diverges almost surely almost everywhere and is  almost surely not a Fourier-Stieltjes series
(cf. \cite{Kahane1985}, p. 54). Regular and irregular properties of the function defined by the series (\ref{eq:RTS}) are studied under the condition $\sum \mathbb{E}(X_n^2) <+\infty $
(cf. \cite{Kahane1985}, \cite{MP1981}).

We would like to study the series (\ref{eq:RTS}) under the   following assumptions 
\begin{equation}\label{eq:assumption}
 \sum_{n=1}^\infty \mathbb{E}X_n^2\  = \infty, \qquad \sum_{n=1}^\infty \mathbb{E}X_n^4\  < \infty. 
 \end{equation} 
One of our objects of study is  the behaviors of the partial sums of the pseudo function  defined by the series 
(\ref{eq:RTS}):
\begin{equation}\label{eq:PS}
S_N(t) : = S_N(t, \omega) := \sum_{n=1}^N X_n \cos  (nt +\Phi_n).
 \end{equation} 
We also assume that $\{\Phi_n\}$ are independent and uniformly distributed  on the circle $\mathbb{T}:=\mathbb{R}/2\pi \mathbb{R}$ which is identified with the interval  $[0, 2\pi)$ and that 
$\{X_n\}$ are independent quasi-gaussian and independent of $\{\Phi_n\}$.  A real random variable $X$ is said to be quasi-gaussian if 
$\mathbb{E} (X^{2m}) = O(K^m m!)$ for some $K\ge 0$. Recall that $X$ is said to be 
subgaussian if $\mathbb{E} e^{\lambda X}\le e^{\frac{1}{2} \tau^2 \lambda^2}$
for some $\tau >0$ and all $\lambda\in \mathbb{R}$. We know that  $X$ is subgaussian if and only if it is centered and quasi-gaussian (\cite{Kahane1985}, p. 82). 
 As we shall prove, the behavior of $S_N(t)$ is very multifractal, meaning that for uncountably many  functions
 $L: \mathbb{N}\to \mathbb{R}$ of different orders, the sets $\{t: S_N(t) \sim L(N)\}$ have positive Hausdorff dimensions.
 
 The partial sum $S_N$ is a stationary process on $\mathbb{T}$ with its correlation function equal to
 $
     \mathbb{E} S_N(t) S_N(s) = H_N(t-s),
 $
 where
 $$  
     H_N(t) = \frac{1}{2}\sum_{n=1}^N  \mathbb{E}X_n^2 \cos  nt.
 $$
 The limit $H(t):= \lim_{N\to \infty} H_N(t)$, if exists, is defined to be the correlation function of the series (\ref{eq:RTS}). 
 In many cases, the correlation function
 $$
     H(t)= \frac{1}{2}\sum_{n=1}^\infty \mathbb{E}X_n^2 \cos  nt
 $$
 is well defined, integrable on $\mathbb{T}$ and continuous everywhere except $t=0$ (cf. Section \ref{sec:chaos}).  The correlation function $H$ will play an important role in the study of the series (\ref{eq:RTS}).  
 
 When $\mathbb{E}X_n^2=\frac{\alpha^2}{n}$, we get a special correlation function
  $$
     H_\alpha (t)= \frac{\alpha^2}{2}\sum_{n=1}^\infty \frac{ \cos  nt }{n}= - \frac{\alpha^2}{2}\log \left(2\left|\sin \frac{t}{2}\right|\right).
  $$
  Its exponential is the $\frac{\alpha^2}{2}$-order Riesz kernel
  $$
      e^{H_\alpha(t)} \asymp \frac{1}{\left|\sin \frac{t}{2}\right|^{\frac{\alpha^2}{2}}}.
  $$

An effective tool that we shall use to study the series (\ref{eq:RTS}) is a class of  the multiplicative chaotic measures, which  are formally defined by 
\begin{equation}\label{eq:CM}
\mu_{\omega}= \prod_{n=1}^\infty \exp\left( X_n \cos  (nt +\Phi_n) - \log \mathbb{E}I_0(X_n)\right) dt,
 \end{equation} 
 where $I_0$ is the modified Bessel function of first order:
 $$
     I_0(\alpha) = \int_0^{2\pi} e^{\alpha \sin  x} \frac{dx}{2\pi}    \qquad(\alpha \in \mathbb{R}). 
 $$
 More precisely, $\mu_{\omega}$ is the weak limit of the partial products of the infinite product in (\ref{eq:CM}). These partial products
 form a non-negative martingale for each fixed $t$, which ensures  the existence of the limit.
 We shall find conditions ensuring the non nullity of the measure $\mu_{\omega}$. A simple such condition is
 \begin{equation}\label{eq:EnergyCond}
      \int\int e^{H(t-s)}dtds <\infty
 \end{equation}
 (cf. Proposition \ref{prop:L2}).
 Under the energy condition (\ref{eq:EnergyCond}), we can define a probability measure $\mathcal{Q}$ on the product space $\mathbb{T}\times \Omega$ by 
 the equality
 $$
      \mathbb{E}_{\mathcal{Q}} h(t, \omega) = \mathbb{E} \int_0^{2\pi} h(t, \omega) d\mu_{\omega}(t)
 $$
 holding for all bounded functions $h$. This measure is called Peyri\`ere measure. An important fact is that ``$\mathcal{Q}$-almost surely'' means ``almost surely (with respect to $P$)
 $\mu_{\omega}$-almost everywhere''. Another important fact is that the random variables $\{X_n \cos(nt+\Phi_n)\}$ defined on $\mathbb{T}\times \Omega$ are $\mathcal{Q}$-independent
 (cf. Theorem \ref{peyriere-thm}).
 On the other hand, the measure $\mu_{\omega}$, through the measure $\mathcal{Q}$, well captures the points $t$ for which the series (\ref{eq:RTS}) has a specific property. 
 In the definition of $\mu_{\omega}$ we can replace $\{X_n\}$ by $\{\alpha X_n\}$ to give a measure $\mu_{\omega, \alpha}$. So, we possess many measures as tools. 
 
 For the typical case $ X_n=\frac{1}{\sqrt{n}}$, we shall prove  that 
$\mu_{\omega, \alpha} =0$ almost surely if $|\alpha| >2$ (cf. Theorem \ref{thm:Degen}). However, for $|\alpha| <2$ the Hausdorff dimension of the measure  $\mu_{\omega, \alpha} $ is almost surely equal to
$$
    \dim \mu_{\omega, \alpha}  = 1- \frac{\alpha^2}{4}       \qquad (|\alpha| <2)
$$
(cf. Theorem \ref{thm:dim}). 
The Peyri\`ere measure is denoted $\mathcal{Q}_\alpha$ when 
$|\alpha|<2$.

Let us state the  following two representative results obtained in this paper on the series (\ref{eq:RTS}). See Theorem \ref{thm:ILL0} and Theorem \ref{thm:LD}. 

\begin{thm}\label{thm:main1}
 Let us consider the series (\ref{eq:RTS}) with $X_n = \frac{1}{\sqrt{n}}$. Let $\alpha \in (-2, 2)$. 
Almost surely  $\mu_{\omega, \alpha}$-almost everywhere we have 
$$
    \lim_{N\to \infty} \frac{1}{\log N} \sum_{n=1}^N \frac{\cos (nt +\Phi_n)}{\sqrt{n}}= \frac{\alpha}{2}.
$$  
Moreover, for any $\eta >0$,  we have the following large deviation
$$
    \lim_{N\to \infty} \frac{1}{\log N} \log \mathcal{Q}_\alpha\left\{ (t, \omega):  \frac{1}{\log N}\sum_{n=1}^N \frac{\cos (nt+\Phi_n)}{\sqrt{n}} \not\in \frac{\alpha}{2}+ [-\eta, \eta]\right\}
    =-\eta^2.
$$
\end{thm}

On the torus $\mathbb{T}^d$, a standard trigonometric multiplicative chaotic measure is defined by 
\begin{equation*}\label{eq:TMC}
 Q_\alpha \sigma =\prod_{n}\exp\left[ \alpha X_n \cos \big(n\cdot x +\Phi_n\big) - \log \mathbb{E} I_0(\alpha X_n)\right]d\sigma(x)
\end{equation*} 
where $\alpha$ is a real parameter, $\sigma$ is a finite Borel measure on the torus $\mathbb{T}^d$, $\{\Phi_n\}$ are iid and uniformly distributed on the torus 
$\mathbb{T}^d$, $\mathbb{E}X_n^2 =\frac{1}{|n|^{d}}$ and the product is taken over $n\in \mathbb{Z}^d\setminus\{0\}$ but only one of $n$ and $-n$ is taken into account.  
\medskip

Let us quote the following result in the case $X_n = \frac{1}{|n|^{d/2}}$ 
(cf. Theorem \ref{thm:dim} and Theorem \ref{thm:dimd}). 

\begin{thm}\label{thm:main2}
Let $\tau(d) = \frac{\pi^{d/2}}{\Gamma(d/2)}$, which is the half area of the unit sphere in $\mathbb{R}^d$. 
For any unidimensional measure $\sigma$ on $\mathbb{T}^d$ we have\\
\indent {\rm (1)} \  $Q_\alpha\sigma =0$ if $\dim \sigma <\frac{\alpha^2}{4} \tau(d)$;\\
\indent {\rm (2)} \  $\dim Q_\alpha \sigma = \dim \sigma- \frac{\alpha^2}{4} \tau(d)$ if $\dim \sigma > \frac{\alpha^2}{4} \tau(d)$.   
\end{thm}

 Our main effort is to determine   the kernel and the image of the chaotic operator  $\mathbb{E} Q_\alpha$
 for which we get the following result (cf. Theorem \ref{thm:ker} and Theorem \ref{thm:kerd}).
 
\begin{thm}\label{thm:main3}
 Let $\tau(d) = \frac{\pi^{d/2}}{\Gamma(d/2)}$. We have 
  $$
  \mathcal{S}_{\frac{\tau(d) \alpha^2}{4} -} \subset {\rm Ker} \mathbb{E} Q_\alpha \subset \mathcal{S}_{\frac{\tau(d) \alpha^2}{4} +};
  $$
   $$
  \mathcal{R}_{\frac{\tau(d) \alpha^2}{4} +} \subset {\rm Im} \ \mathbb{E} Q_\alpha \subset \mathcal{R}_{\frac{\tau(d) \alpha^2}{4} -}.
  $$
  \end{thm}
 
 The spaces $\mathcal{S}_{a}, \mathcal{R}_{a}$ etc will be discussed in  Section \ref{sec:chaos}. Let us point out that Theorem \ref{thm:main2} is actually a consequence of Theorem \ref{thm:main3}  and the decomposition principle stated in Theorem \ref{decomposition-principle}. It is a very interesting problem to completely determine the image and the kernel of  $\mathbb{E} Q_\alpha$.
 Associated to the percolation problem on a tree, there is a multiplicative chaotic operator
 for which the image and the kernel are completely determined (cf. \cite{Fan1990}. See \cite{Fan1989b} for the detailed proof.) 
 \medskip

Our work will actually concentrate on the study of trigonometric multiplicative chaos. This enters into the general theory of multiplicative chaos of Kahane \cite{Kahane1987a}.
 Other multiplicative chaos  have been already studied, including gaussian chaos \cite{Kahane1985b}, L\'evy stable chaos \cite{Fan1997}, Dvoretzky covering chaos \cite{Kahane1987b}, tree percolation chaos
 \cite{Fan1990}  etc.  
General infinitely divisible chaos, which includes the L\'evy stable chaos,  are studied in \cite{BM2003}.
 There is recently an active study, which produces  a large literature, on gaussian chaos because of its link to physics, especially to  Liouville quantum fields.
We just cite here some of these works \cite{BJRV2013,BKNSW2015,DS2011,GM2021,RSV2014,Shamov2016} and invite the reader to refer the references therein
and the survey papers \cite{RV2014,RV2016}.  We point out that for both Dvoretzky covering chaos and tree percolation chaos and only for them
\cite{Kahane1987b,Fan1990} (see also \cite{Fan1989b}), we know the exact kernel and image 
of the corresponding chaotic projection operator $\mathbb{E}Q$ (see \S 2 for the definition of this operator $\mathbb{E}Q$).
There are also an intensive study of different aspects of particular chaotic measures and their applications \cite{Barral1999,Barral2001,BF2005,Ben1990, Fan1990, Fan1993, Fan1995, Fan1995b,  Fan2002, Fan2002b, Fan2004, FK1993,FK2001, FK2021, Liu2000, LR2000}. 

 \medskip
 
 In the section \S 2, we give a brief recall of the general theory of multiplicative chaos and state all known basic results that we shall use.
 Our trigonometric chaos on $\mathbb{T}$ are constructed in the section \S 3, where the correlation functions and  the associated  kernels are examined,     the capacity  and the dimension of measures are recalled, as well as their relations. The section \S4 is devoted to $L^2$-theory and its consequences. There we also discuss the decomposition and the composition of our typical chaotic operators.
 The degeneracy of chaotic operator is investigated in the section \S 5, while the kernel and the image of the projection
 $\mathbb{E} Q_\alpha$ are studied in the section \S 6.  
 All these can be generalized to the torus $\mathbb{T}^d$ with $d\ge 1$ as we shall show in the section \S 7. 
 Our chaotic measures are used in the section \S 8 to study the divergence of the series (\ref{eq:RTS}).  
  \medskip

 {\bf Notation.} We shall adopt the following notation. Let $u(x)$ and $v(x)$ be two functions. We write $u(x) \ll v(x)$ if there exists a positive constant $C>0$
 such that $u(x) \le C v(x)$ for $x$ in the domain of definition of   $u$ and $v$. If $u(x)\ll v(x)$ and $v(x)\ll u(x)$, we write $u(x) \asymp v(x)$.
 When $\lim_{x\to x_0}\frac{u(x)}{v(x)} =1$, we write $u(x) \sim v(x)$ ($x\to x_0$) or simply  $u(x) \sim v(x)$.

\medskip
 
{\it Acknowledgement and Addendum.}  
We would like to thank J. Barral and V. Vargas for pointing out the reference \cite{Junnila2020} in which Janne Junnila obtained  the result about the full action
on the Lebesgue measure, even for $\alpha$ in a complex domain containing $(-2, 2)$ in the case of Theorem \ref{thm:main1} and he also studied other fields than trigonometric field on $\mathbb{T}$.  
The first author is partially supported by NSFC (grant no.11971192). 

\section{General multiplicative chaos}
In this section, we recall the basic results in the theory of $T$-martingales, or of multiplicative chaos.
The origin goes back to \cite{Kolmogorov1961,Mandelbrot1971,Mandelbrot1972,Mandelbrot1974a,Mandelbrot1974b} where understanding turbulence was the motivation. The general theory was developed by Kahane \cite{Kahane1987a}.
First seminal works on the subject were \cite{KP1976} and \cite{Kahane1985b}. 
In the next section, 
we shall introduce our $\mathbb{T}$-martingales, that we shall study in this paper.
Let us point out that the theory in \cite{Kahane1987a} is  generalized in \cite{WW1996} to a case  including some  dependent multiplicative cascades. 
The materials presented below come from \cite{Kahane1987a, Fan1991, FK2008}. 

\subsection{Construction of $T$-martingales}
Let $(T, d)$ be a compact (or locally compact) metric space and
$(\Omega, \mathcal{A}, P)$
be a probability space. We are given an increasing sequence
$\{\mathcal{A}_n\}_{n \ge 1}$ of sub-$\sigma$-fields of
$\mathcal{A}$ and a sequence of random functions
$\{P_n(t, \omega)\}_{n \ge 1}$ ($t \in T, \omega \in \Omega$). 
We make the following assumptions:
\begin{itemize}
  \item[(H1)]  \ \ \ $P_n(t)=P_n(t, \omega)$ are non-negative and independent
  processes;
     $P_n(\cdot, \omega)$ is Borel measurable for almost all $\omega$;
     $P_n(t, \cdot)$ is $\mathcal{A}_n$-measurable for each $t$.
  \item[(H2)] \ \ \
    $\mathbb{E} P_n(t) =1$ for all $t \in T$.
\end{itemize}
Such a sequence $\{P_n\}$ is called a sequence of {\em weights}
adapted to $\{\mathcal{A}_n\}$.

 Let
\begin{equation}\label{eq:QM}
     Q_n(t) = Q_n(t, \omega)= \prod_{j=1}^n P_j(t, \omega).
\end{equation}
For any $t\in T$, $\{Q_n(t)\}$ is a martingale. We call $\{Q_n(t)\}$ ($t\in T$) a $T$-martingale, or a martingale indexed by $T$.
For any $n\ge 1$ and any positive Radon measure $\sigma$ on $T$
(we write $\sigma \in \mathcal{M}^+(T)$), we consider the random
measures $Q_n\sigma$ defined by
$$
      Q_n\sigma (A) = \int_A Q_n(t) d \sigma(t)
      \qquad (A \in \mathcal{B}(T))
$$
where $\mathcal{B}(T)$ is the Borel field of $T$.
It is clear that for any $A \in \mathcal{B}(T)$, $Q_n\sigma (A)$ is a positive
martingale, so it converges almost surely (a.s. for short).
So does $\int \phi(t) d Q_n \sigma (t)$ for any bounded Borel function $\phi$.

\subsection{Basic results in the theory of multiplicative chaos}
The following fundamental theorem is proved based on  the last fact stated above.

\begin{thm} {\rm (\cite{Kahane1987a})}
    For any Radon measure $\sigma \in \mathcal{M}^+(T)$,
    almost surely
    the random measures $Q_n\sigma$ converge weakly  to
    a random measure $S$, which will be denoted $Q\sigma$.
\end{thm}

We may consider $Q$ as an operator which maps measures into
random measures.
We call $Q$ a {\em multiplicative chaotic operator} or simply {\em chaotic operator}, and
$Q\sigma$ a {\em multiplicative chaotic measure} or simply {\em chaotic measure} which, in some special cases,
describes the limit energy state of a turbulence
\cite{Mandelbrot1974a,Mandelbrot1974b,Kahane1985,Fan1989a}.


There are two possible extreme cases. The first one is that
$Q\sigma=0$ a.s.
(the energy is totally dissipated). The second one is that 
$Q_n \sigma(T)$ converges in $L^1$ or equivalently  $\mathbb{E} Q\sigma (T)= \sigma(T)$
(the energy is conserved).
If the first case occurs, we say that $Q$ {\em degenerates} on $\sigma$
or $\sigma$ is  $Q$-{\em singular}.
If the second case occurs, we say that $Q$ {\em fully acts} on $\sigma$
or $\sigma$ is  $Q$-{\em regular}.

We define a map $\mathbb{E}Q: \mathcal{M}^+(T)\to \mathcal{M}^+(T)$
by
$$
    (\mathbb{E}Q\sigma)(A) = \mathbb{E} (Q\sigma(A))
    \qquad (A \in \mathcal{B}(T)).
$$
That $\sigma$ is  $Q$-singular (resp. $Q$-regular) is equivalent to
$\mathbb{E}Q \sigma =0$ (resp.
$\mathbb{E}Q \sigma =\sigma$).

\begin{thm} {\rm (\cite{Kahane1987a})}\label{Decomposition}
    Any Radon measure $\sigma \in \mathcal{M}^+(T)$
    can be uniquely decomposed into
    $\sigma= \sigma_r + \sigma_s$ where $\sigma_r$ is a $Q$-regular measure
    and $\sigma_s$ is a $Q$-singular measure.
    Both $ \sigma_r $ and $\sigma_s$ are restrictions of $\sigma$, that is to say
    $\sigma_r= \sigma 1_B$ for some Borel set $B$, so that $\sigma_s=\sigma 1_{B^c}$.
\end{thm}

The operator $\mathbb{E}Q$ extended to the space  $\mathcal{M(T)}$
 is thus a projection whose image (resp. kernel)
    consists of $Q$-regular (resp. $Q$-singular) measures.
We call $\mathbb{E}Q$ the chaotic projection operator.

We are concerned with  properties of the random measure $Q\sigma$, of the
 operator $Q$ or of the projection $\mathbb{E}Q$. Here are some fundamental questions:
\begin{itemize}
\item[] {\bf Question 1}  Does $Q$ degenerate on $\sigma$?
\item[] {\bf Question 2} Does $Q$ fully act on $\sigma$?
\item[] {\bf Question 3} What is the dimension of the measure $Q\sigma$?
\item[] {\bf Question 4} What are the possible relations between two measures $Q'\sigma'$ and $Q''\sigma''$
for two different operators $Q'$ and $Q''$ defined in the same way as $Q$? For example, when $Q'\sigma'$ and $Q''\sigma''$
are mutually singular or absolutely continuous one with respect to the other ?
\end{itemize}
In the following, we state some results in the general case. Either they are
partial answers
to one of these questions, or they provide some useful tools. 

\begin{thm}[\cite{Kahane1987a}] \label{thm:Kdegen}
   Suppose that $\mathcal{H}^\alpha (T) <\infty $
   where $\mathcal{H}^\alpha$ denotes the $\alpha$-dimensional
   Hausdorff measure and that there exist constants $0<h<1$
   and
   $C>0$ with the property that for any ball $B$ with radius $r$
   there exists an integer $n=n(B)$ such that
   \begin{equation}~\label{singularity}
       \mathbb{E} \left(\sup_{t \in B} Q_n(t)\right)^h
       \le C r^{\alpha(1-h)}.
   \end{equation}
   Then all Radon measures on $T$ are $Q$-singular.
\end{thm}

This provides a good tool to verify the $Q$-singularity of $\sigma$.
In fact, the condition $\dim \sigma <\alpha$
together with (\ref{singularity}) implies the $Q$-singularity of $\sigma$.

On the other hand,
the following gives a simple condition of $Q$-regularity, which is the condition for the $L^2$-convergence of the martingale $Q_n\sigma(\mathbb{T})$.

\begin{thm}[\cite{Kahane1987a}]
  The operator $Q$ fully acts on $\sigma$ and  $\mathbb{E} (Q\sigma(T))^2 <\infty$
  if and only if
   \begin{equation}~\label{L2regularity}
     \lim_{N\to \infty}  \int \int \prod_{n=1}^N \mathbb{E} P_n(t)P_n(s) d \sigma(t) d \sigma(s)
       <\infty.
   \end{equation}
\end{thm}


Suppose that $\sigma$ is a $Q$-regular probability measure.
A useful tool for studying the measure $Q\sigma$ is the {\em Peyri\`ere measure}
$\mathcal{Q}$  on the product space $T\times \Omega$ defined by
\begin{equation}~\label{peyrieremeasure}
      \int_{T\times \Omega} \varphi(t, \omega) d \mathcal{Q}(t, \omega)
          = \mathbb{E} \int_T \varphi(t, \omega) d  Q\sigma(t)
\end{equation}
for non-negative measurable functions $\varphi$.
\medskip

\begin{thm}[\cite{Kahane1987a}]~\label{peyriere-thm}
      Suppose that $\sigma$ is a $Q$-regular probability measure
      and that the probability law of the variable $P_n(t)$ 
      is independent of $t$. Then $P_n$'s, considered as random
      variables with respect to $\mathcal{Q}$, are $\mathcal{Q}$-independent.
\end{thm}

If the conditions in Theorem \ref{peyriere-thm} are satisfied, 
for any bounded function or positive function $\Phi$, we have
\begin{equation}\label{eq:ExpectQ}
   \mathbb{E}_\mathcal{Q} \Phi(P_n(t)) = \mathbb{E} \Phi(P_n(t)) P_n(t),
\end{equation}
where the term on the right hand side is independent of $t$.
\medskip

A direct application of the Peyri\`ere measure leads to that
almost surely for $Q\sigma$-almost every $t \in T$ we have
\begin{equation}
         \lim_{n \to \infty} \frac{1}{n}
         \sum_{k=1}^n \log P_k(t, \omega)
            = \lim_{n \to \infty} \frac{1}{n}
         \sum_{k=1}^n  \mathbb{E} P_k\log P_k.
\end{equation}
This result
may be used to study the dimension and the multifractality of the measure
$Q\sigma$.
\medskip

Now suppose that we are given two sequences of
weights $\{P_n'(t)\}_{n \ge 1}$ and $\{P_n''(t)\}_{n \ge 1}$
 adapted to $\{\mathcal{A}_n'\}$ and
$\{\mathcal{A}_n''\}$, respectively defined on
probability spaces $(\Omega',\mathcal{A}', P')$
and $(\Omega'',\mathcal{A}'', P'')$.
It is obvious that $\{P_n\}$
defined by $P_n(t) = P_n''(t) P_n'(t)$
is a sequence of weights adapted to
$\{\mathcal{A}_n' \otimes\mathcal{A}_n''\}$, defined on the
product space $(\Omega, \mathcal{A}, P) = (\Omega' \times \Omega'',
\mathcal{A}'\otimes \mathcal{A}'', P'\otimes P'')$.
We denote by $Q', Q''$ and $Q$ the three operators corresponding
to the above three sequences of weights.
The following decomposition principle establishes a relationship between
$Q', Q''$ and $Q$.
\medskip

\begin{thm}[\cite{FK2008}]~\label{decomposition-principle}
      Under the above condition,  we have
\begin{itemize}
  \item [(a)] \mbox{\rm a.s.} $Q\sigma = Q'' (Q'\sigma)  
              $ for any measure
      $\sigma \in \mathcal{M}^+(T)$.
  \item [(b)] $ \sigma \in \mbox{\rm Im} \mathbb{E}Q \Rightarrow
           Q'\sigma \in \mbox{\rm Im} \mathbb{E}Q''$ for
           almost all $\omega' \in \Omega'$.
  \item [(c)]
         $ \sigma \in \mbox{\rm Ker} \mathbb{E}Q \Rightarrow
           Q'\sigma \in \mbox{\rm Ker} \mathbb{E}Q''$
           for almost all $\omega' \in \Omega'$.
   \item [(d)] $\mathbb{E} Q = \mathbb{E}Q'' \mathbb{E}Q'$.
\end{itemize}
\end{thm}

\medskip

Let $Q'$ and $Q''$ be two operators associated respectively  to
$\{ P_n'\}$ and $\{ P''_n \}$. Now we do not suppose the independence of
$\{P_n'\}$ and $\{P''_n\}$. But we suppose that the law
of the vector $(P'_n(t),P''_n(t))$ is independent of $t$.
We have a Kakutani type criterion for the dichotomy of the mutual
absolute continuity of $Q'\sigma$ and $Q''\sigma$.
\medskip

\begin{thm}[\cite{Fan1991}]~\label{ThmFan1991} Assume the above assumptions. Suppose furthermore that
$\sigma$ is a $Q'$-regular probability measure.  We have
\begin{itemize}
  \item [(a)] $\prod_{n=1}^\infty \mathbb{E} \sqrt{P'_n P''_n} >0 \Rightarrow \ a.s.\ 
      Q'' \sigma \ll Q'\sigma$ and $\sigma$ is $Q''$-regular.
  \item [(b)] $\prod_{n=1}^\infty \mathbb{E} \sqrt{P'_n P''_n} =0 \Rightarrow \ a.s.\ 
      Q'' \sigma \perp Q'\sigma$.
\end{itemize}
\end{thm}

\section{Trigonometric multiplicative chaos on $\mathbb{T}$}
\label{sec:chaos}
Now, in this section, we  construct a class of $\mathbb{T}$-martingales, which define our trigonometric multiplicative chaos.
We also discuss the associated correlation functions and their relations to potential theory and dimension theory.

\subsection{Definition of $\mathbb{T}$-martingales}

We always make the following assumptions:
\begin{itemize} 
\item [(H1)]  $X_n =\alpha_n Y_n$ where $Y_n$'s are normalized independent  quasi-gaussian real random variables. That means
\begin{equation*}\label{eq:CondOnY}
   \mathbb{E} Y_n^2=1, \quad \mathbb{E} Y_n^{2m} = O(K^m m!) 
 \end{equation*}
 for  some  $K>0$ and  for  all $m \ge 2$ and $n\ge 1$;
 \item[(H2)] $\Phi_n$'s, which will be denoted $\omega_n$,  are independent random variables which are uniformly distributed on $\mathbb{T}$;
 \item [(H3)] the real coefficients $\alpha_n$'s satisfy
  \begin{equation*}\label{eq:TM-assumptions}
    \sum_{n=1}^\infty |\alpha_n|^2=\infty, \qquad   \sum_{n=1}^\infty |\alpha_n|^4<\infty.
    \end{equation*}
 \end{itemize}
 We define the weights
\begin{equation}\label{eq:W1}
   P_n(t) = \frac{e^{\alpha_n Y_n \cos  (nt+\omega_n)}}{\mathbb{E}e^{\alpha_n Y_n \cos  (nt+\omega_n)}}= \frac{e^{\alpha_n Y_n \cos  (nt+\omega_n})}{\mathbb{E}I_0(\alpha_n Y_n)},
\end{equation}
where $I_0$ denotes the modified Bessel function of first kind.
Then we define  the $\mathbb{T}$-martingale $\{Q_n(t)\}$ by 
 \begin{equation}\label{eq:TM}
   Q_n(t) = \exp \left[\sum_{k=1}^n \Big(\alpha_k Y_k\cos  (nt +\omega_k)- \log \mathbb{E} I_0(\alpha_k Y_k)\Big)\right].
 \end{equation}
 We shall denote by $Q$ the multiplicative operator defined by this  $\mathbb{T}$-martingales $\{Q_n(t)\}$.

For  this  multiplicative chaotic operator $Q$, 
we define its correlation function by 
  \begin{equation}\label{eq:H}
       H(t) = \frac{1}{2}\sum_{n=1}^\infty \alpha_n^2 \cos  nt,
     \end{equation}
    and its associated kernel
     \begin{equation}\label{eq:Fi}
        \Phi(t) = e^{H(t)}.
     \end{equation}

    The following simple properties of the Bessel function $I_0$ will be frequently used in the sequel. 
    
  \begin{lem}\label{lem:Bessel}    Suppose that $X$ is quasi-gaussian and $\mathbb{E} X^2 =1$. 
  As $\alpha \to 0$ we have
  \begin{eqnarray*}
     \mathbb{E}I_0\left(\alpha X\right) & = &  1 + \frac{\alpha^2}{4} + O(\alpha^{4}) = e^{\frac{\alpha^2}{4} + O(\alpha^{4})}, \\
       \frac{d}{d\alpha} \mathbb{E}I_0\left(\alpha X\right) & = &  \frac{\alpha}{2}   + O(\alpha^{3}), \\
         \frac{d^2}{d\alpha^2} \mathbb{E}I_0\left(\alpha X\right) & = &  \frac{1}{2} +  O(\alpha^2).
  \end{eqnarray*}
  \end{lem} 
  \begin{proof}
   From the well known fact
  $
       I_0 (x)  =\sum_{m=0}^\infty \frac{1}{m!^2} \left(\frac{x}{2}\right)^{2m} 
  $,  we get
  $$
      \mathbb{E}I_0\left(\alpha X\right)  = \sum_{m=0}^\infty \frac{1}{m!^2}  \frac{\alpha^{2m}}{2^{2m}}  \mathbb{E} (X^{2m})  = 1+ \frac{\alpha^2}{4} + \sum_{m=2}^\infty \frac{1}{m!^2}  \frac{ \alpha^{2m} }{2^{2m}} 
      \mathbb{E} (X^{2m}) .
  $$
  Since $X$ is quasi-gaussian, we have $ \mathbb{E} (X^{2m})  \le M K^m m!$ for some constants $M>0$ and $K>0$. 
  It follows that the last sum is bounded by $ \frac{M K^2}{4} \alpha^4 e^{\frac{K}{4} \alpha^2}$. The two other estimates can be similarly obtained. 
  \end{proof}
  
  The next lemma states a basic relation between the chaotic operator and its correlation function.

  \begin{lem}\label{lem:EQQ} Under the assumptions  (H1), (H2) and (H3) 
  we have 
  $$
     \mathbb{E} [Q_n(t)Q_n(s)]  \asymp \exp \left( \frac{1}{2}\sum_{k=1}^n \alpha_k^2 \cos  k (t-s) \right). 
     $$
  \end{lem}
  \begin{proof} For any real number $\lambda$ we have the equality
  $$
     \mathbb{E} \exp\left(\lambda [\cos (k t +\omega_k) + \cos (k s + \omega_k)]\right) =  I_0(2\lambda \cos \pi k(t-s)).
  $$
  Indeed, 
  using the formula $\cos \alpha + \cos \beta = 2 \cos \frac{\alpha +\beta}{2} \cos \frac{\alpha -\beta}{2}$, we get 
  $$
    \cos (k t +\omega_k) + \cos (k s + \omega_k) = 2  \cos  (k(t+s)/2 + \omega_k) \cos \pi k(t-s).
  $$
  Then, the equality follows from the definition of $I_0$ and the translation invariance of Lebesgue measure.
 It follows that   
  \begin{eqnarray*}
    \mathbb{E} [Q_n(t) Q_n(s)] 
    &=& \prod_{k=1}^n \frac{\mathbb{E}I_0(2\alpha_kY_k\cos \pi k(t-s))}{[\mathbb{E}I_0(\alpha_kY_k)]^2}.
     \end{eqnarray*}
 But, according to Lemma \ref{lem:Bessel},  we have  
 \begin{eqnarray*}     
 \frac{\mathbb{E}I_0(2\alpha_k Y_k\cos \pi k (t-s))}{[\mathbb{E}I_0(\alpha_k Y_k)]^2}
    & =& \exp \left[  \alpha_k^2\Big(\cos^2 \pi k (t-s) - \frac{1}{2}\Big)  + O(\alpha_k^{4}) \right].
 \end{eqnarray*}
 Now we can conclude because $\cos^2 x -\frac{1}{2} = \frac{1}{2}\cos 2x$ and $\sum \alpha_k^4 <\infty$.
\end{proof}    

Lemma \ref{lem:EQQ} is fundamental for all computations  in the sequel. Remark that the estimate in Lemma \ref{lem:EQQ}
only depends on $\alpha_k^2$ but not on the distributions of $Y_k$. So, without loss of generality we can assume that $Y_k =1$ for all $k$. In other words, 
we can treat all sequences $\{Y_k\}$ satisfying the assumption (H1) just as in the special case of $Y_k=1$. 
From now on, we assume that $Y_k=1$.

\medskip

A typical case is $\alpha_n = \frac{\alpha}{\sqrt{n}}$, for which the correlation function is equal to 
   $$
         H_\alpha(t) = - \frac{\alpha^2}{2} \log |\sin \pi t| +O(1).
   $$

We finish our construction by pointing out that instead of  $\cos x$ we can consider any $2\pi$-periodic function
     $f(x)$ such that $$J(\alpha):=\int_0^{2\pi} e^{\alpha f(x)} \frac{dx}{2\pi} <\infty \quad {\rm for} \ \ \alpha \in (-\delta, \delta)$$ 
     for some $\delta>0$. This function $J$ plays the role of the Bessel function $I_0$.  Let $m_n =\int_0^{2\pi} f(x)^m \frac{dx}{2\pi}$ be the moments ($n\ge 0$). We have 
     $$
           J(\alpha) = \sum_{n=0}^\infty \frac{m_n}{n!} \alpha^n.
     $$
     Some conditions on the moments of $f$ are needed. 

\subsection{Correlation function $H$ and Kernel $\Phi$}
Here we first present some conditions for $H$ to be pointwise well defined and for $H$ even $\Phi$ to be integrable. The function
    $\Phi$ will play the role of kernel in the sense of potential theory.   We refer to \cite{Zygmund}, Vol. 1, Chapter V.
 
    It is well known that if $\alpha_n^2\downarrow 0$, then the cosine  series (\ref{eq:H}) defining $H$ converges uniformly on every
    interval $[\delta, 2\pi-\delta]$ ($\delta>0$) and consequently the function $H$ is continuous in the open interval  $(0,2\pi)$.

    A sequence $\{a_n\}$ is  said to be convex  if $a_{n+2}-2a_{n+1}+a_n\ge 0$  for all $n$. 
    A convex sequence tending  to $0$ must be decreasing. It is also well known that
     if $\{\alpha_n^2\}$ is convex and $\alpha_n^2 \downarrow 0$, then the function $H$ defined by the series (\ref{eq:H}) is lower bounded and Lebesgue-integrable. 
     In the following lemma, we give a condition ensuring that the partial sums of  (\ref{eq:H}) is bounded by the sum of  (\ref{eq:H}).

     \begin{lem}\label{lem:upperbound}
      Suppose  $\{\alpha_n^2\}$ is convex and 
      \begin{equation}\label{eq:bdcond}
      \alpha_n^2 \downarrow 0,\quad \alpha_n^2 =O(n^{-1}), \quad
      \alpha_n^2 - \alpha_{n+1}^2 = O(n^{-2}).
      \end{equation}
    Then the partial sums of the series (\ref{eq:H}) are uniformly upper bounded by $H(x)+C$:
    $$
        \forall N\ge 1, \forall x\in \mathbb{T}, 
        \qquad   \frac{1}{2}\sum_{n=1}^N \alpha_n^2 \cos  n t \le H(t) +C
    $$
    where $C$ is a constant. 
    \end{lem}
    \begin{proof}  We repeat the classical proof (cf \cite{Bary1964}, p. 92) and make the observation for the boundedness.
    For simplicity, let $a_n =\frac{\alpha_n^2}{2}$. Choose a positive $a_0$ such that $\{a_n\}_{n\ge 0}$ is convex. Let
    $$
              S_N (x) =\frac{a_0}{2} +\sum_{n=1}^N a_n \cos  n t. 
    $$ 
    By a double application of Abel's transformation, we get
    $$
          S_N(x) = \sum_{n=0}^{N-2} (n+1) \Delta^2a_n K_n(x) + N \Delta a_{N-1} K_{N-1}(x) + a_N D_N(x),
    $$
    where $D_n$ is the Dirichlet kernel and $K_n$ is the F\'ejer kernel, both of them are uniformly bounded by $n$.  
    By the assumption (\ref{eq:bdcond}), we have
    $$
       |N \Delta a_{N-1} K_{N-1}(x) + a_N D_N(x)|\le C
    $$
    for some constant $C>0$. Now we can conclude, because $K_n$ are non-negative and 
    $$
       \sum_{n=0}^\infty (n+1) \Delta^2a_n K_n(x) - \frac{a_0}{2} = H(t).
    $$
    \end{proof}

    The conditions in Lemma \ref{lem:upperbound} are satisfied by 
    $$
        \alpha_n^2 = \frac{1}{n \log^\beta n} \qquad (n \ge 3, \beta \ge 0).
    $$
    ($\alpha_1^2, \alpha_2^2$ are conveniently chosen). 
    
      If $\alpha_n^2 \downarrow 0$, the function $H$ is continuous for $x\not=0$. Also observe that $H$ is even. So, in this case,
      for any probability Borel measure $\sigma$ on $\mathbb{T}$ we have the equivalence
    $$
       \int\int e^{H(t-s)} d\sigma(t) d\sigma(s)<\infty \Longleftrightarrow \int \int_{|t-s|\le \delta} e^{H(t-s)} d\sigma(t) d\sigma(s)<\infty.
    $$ 
    for some or all $\delta>0$.
    If $\sigma$ is the Lebesgue measure, after a change of variables we get 
     $$
       \int\int e^{H(t-s)} dt ds <\infty  \Longleftrightarrow \int e^{H(t)} dt <\infty \Longleftrightarrow \int_0^\delta e^{H(t)} dt <\infty.
    $$ 
    for some or all $\delta>0$. 
    
    Thus, in order to get the integrability of  the kernel $\Phi$,  we only need to known the behavior of $H$ at $t=0$. 
    Let us quote the following result for a family of functions $H$.
    Recall that a positive function $b(u)$ for $u\ge u_0$ is said to be slowly varying if for any $\delta>0$, $b(u) u^\delta$
    is increasing and $b(u) u^{-\delta}$
    is decreasing for $u$ large enough. For any real number $\beta$, the function $\log^\beta u$ ($u\ge 3$) is slowly varying. 
     
      \begin{lem}[\cite{Zygmund}, p.187-188]\label{lem:bound-at-0} Let $f_\tau(t) = \sum_{n=1}^\infty \frac{b_n}{n^\tau} \cos nt$ where $0<\tau\le1$
      and $b_n =b(n)$ for some slowly varying function $b(\cdot)$. When $0<\tau<1$, we have
      $$
           f_\tau(t) \sim \Gamma(1-\tau) \sin \frac{\pi \tau}{2}\cdot \frac{b(t^{-1})}{t^{1-\tau}}, \quad t \to +0.
      $$
      When $\tau =1$ and $\sum \frac{b_n}{n}=\infty$, we have
      $$
            f_1(t) \sim \int_1^{t^{-1}}\frac{b(u)}{u}du,   \quad t \to +0.
      $$
      \end{lem}

      Let us consider the functions 
   $$
      h_{\alpha, \tau, \beta}(t) = \sum_{n=3}^\infty \frac{\alpha}{n^\tau \log^\beta n} \cos  nt \quad (\alpha>0, 0<\tau \le 1, \beta\in \mathbb{R}).
   $$
   
   \begin{cor}
   If $0<\tau<1$, we have
   $$
    h_{\alpha, \tau, \beta}(t) \sim \frac{C_{\alpha, \tau}}{t^{1-\tau} \log^\beta\frac{1}{t}}, \quad t \to +0
    $$
   with $C_{\alpha, \tau} =\alpha \Gamma(1-\tau) \sin\frac{\tau \pi}{2}$,  and $e^{h_{\alpha, \tau, \beta}(t)}$ is not integrable on $\mathbb{T}$.
    \end{cor}
    
     \begin{cor}
    If $\tau=1$ and $0<\beta<1$, we have
   $$
          h_{\alpha, 1, \beta} (t) \sim \alpha \int_3^{t^{-1}} \frac{d u}{u \log^\beta u} \sim \frac{\alpha}{1-\beta} \log^{1-\beta} \frac{1}{t}, \quad t \to +0
   $$
   and $e^{\alpha h_{1, \beta}(t)}$ is integrable on $\mathbb{T}$.  If $\tau=1$ and $\beta>1$, $h_{\alpha, 1,\beta}$ is bounded.
     \end{cor}
     
    \begin{cor}
    If $\tau=1$ and $\beta = 1$, we have
    $$
          h_{\alpha, 1,1} (t) \sim \alpha \int_3^{t^{-1}} \frac{d u}{u \log u} \sim \alpha \log \log  \frac{1}{t}, \quad t \to +0
   $$
    and $e^{\alpha h_{\alpha, 1, 1}(t)}$ is integrable on $\mathbb{T}$.
      \end{cor}
     
        \begin{cor}
     If $\tau=1$ and $\beta  =0$, we have
    $$
          h_{\alpha, 1,0} (t) \sim \alpha \int_3^{t^{-1}} \frac{d u}{u} \sim  \alpha \log  \frac{1}{t}, \quad t \to +0
   $$
    and $e^{\alpha h_{\alpha, 1, 0}(t)}$ is integrable on $\mathbb{T}$ if and only if  $\alpha<1$.
      \end{cor}
      
   \subsection{Capacity, Dimension and Kahane decomposition} \label{sect:capacity}
   We present here some notions from potential theory and their relations to the dimensions of measures.
  Consider the Riesz kernel
  $$
      \Phi_\alpha(t) = e^{H_{\alpha}(t)}  \asymp \frac{1}{ \left|\sin \frac{ x}{2}\right|^\alpha} .
  $$
  The $\alpha$-energy of a finite Borel  measure $\sigma$ is defined by
  $$
        I_\alpha^\sigma =  \int\int \frac{d\sigma(t) d\sigma(s)}{|\sin \frac{t-s}{2}|^{\alpha}}.
  $$
We have forgotten the constant $\frac{1}{2^\alpha}$, because it is usually the finiteness of the energy which plays the role. 

Generally, for a given kernel $\Phi$, 
the $\Phi$-energy of $\sigma$ is similarly defined by
 $$
        I_\Phi^\sigma =  \int\int \Phi(t-s) d\sigma(t) d\sigma(s)
  $$
  and the $\Phi$-potential of $\sigma$ is defined by 
  $$
      U_\Phi^\sigma(t) = \int \Phi(t-s) d\sigma(s).
  $$
  If $\Phi \in L^1(\mathbb{T})$ is non-negative, even and  convex on $(0, 2\pi)$, the $\Phi$-potential theory is well developed (cf.  \cite{KS1963}).  
  The Fourier coefficients $\widehat{\Phi}(n)$ are non-negative and we have the following formula for the energy in terms  of Fourier coefficients: 
  $$
      I^\sigma_\Phi = 4\pi^2 \sum_{n=-\infty}^\infty \widehat{\Phi}(n) |\widehat{\sigma}(n)|^2.
  $$
 The $\Phi$-capacity of a set $E \subset \mathbb{T}$, which is an important notion,  is defined by 
 $$
 {\rm cap}_\Phi(E) = \frac{1}{\inf I^{\sigma}_\Phi},$$
  where the infimum is taken over all probability measures concentrated in $E$. 
  When $\Phi=\Phi_\alpha$, we shall write $I_\alpha^\sigma$ for $I_\Phi^\sigma$.  $U_\alpha^\sigma$ and ${\rm cap}_\alpha(E)$
are similarly understood.   
\medskip
  
 The potential $U_\alpha^\sigma(x)$ describes well the local behavior of the measure $\sigma$ at $x$:
 \begin{equation}\label{eq:U-Lip}
        A \sup_{r>0} \frac{\sigma(B(x, r))}{r^\alpha} \le U_\alpha^\sigma(x) \le B \sup_{r>0} \frac{\sigma(B(x, r))}{r^{\alpha+\epsilon}} 
 \end{equation}
 where $B(x, r)$ denotes the ball centered at $x$ of radius $r$ (an interval of length $2r$), $A$ and $B$ are constants depending only on $\alpha$ and $\epsilon>0$
 (cf. \cite{BFP2010}, p. 103). So, it is easy to see that
  when $I_{\alpha}^\sigma <\infty$, we have $\dim_*  \sigma \ge \alpha$ where $\dim_* \sigma$ denotes the lower Hausdorff dimension of the measure $\sigma$: 
 $$    
    \dim_* \sigma = \inf \{ \dim E: \sigma(E)>0\},
 $$
 where $\dim E$ denotes the Hausdorff dimension of the set $E$ (cf. \cite{Fan1994} for definition).  In \cite{Fan1994},  the upper Hausdorff dimension of
 the probability measure $\sigma$ is defined by
 $$
     \dim^* \sigma = \inf \{ \dim E: \sigma(E)=1\}.
 $$
 If $\dim_* \sigma=\dim^* \sigma$, we write $\dim \sigma$ for the common value, called the Hausdorff dimension of $\sigma$.
 In this case, we say $\sigma$ is unidimensional. 
 Also recall that the lower local dimension of $\sigma$ at $x$ is defined by 
 $$
     \underline{D}(\sigma, x) = \liminf_{r\to 0} \frac{\log \sigma(B(x, r))}{\log r}. 
 $$
 Here we can discretize $r$ by taking a sequence $\{r_n\}$ such that $r_n \downarrow 0$ and $\log r_n \sim \log r_{n+1}$.
 It is proved that \cite{Fan1994}
 $$
      \dim_*\sigma = {\rm ess}\, \inf \underline{D}(\sigma, x), \qquad   \dim^*\sigma = {\rm ess}\, \sup \underline{D}(\sigma, x).
 $$

A finite Borel  measure $\sigma$ is said to be  $\alpha$-regular, if $\sigma$ is a countable sum $\sum \sigma_i$ (convergent in the norm of total variation
and $\sigma_i$ having disjoint Borel supports)
such that $I^{\sigma_i}_\alpha <\infty$ for all $i$; it is said to be  $\alpha$-singular if $\sigma$ is supported by a set of zero $\alpha$-capacity. 
We denote by $\mathcal{R}_\alpha$ (resp. $\mathcal{S}_\alpha$) the set of all $\alpha$-regular (resp. $\alpha$-singular) measures. 
We have the relations
$$
   \sigma\in \mathcal{S}_\alpha \Longrightarrow \dim^*\sigma \le \alpha; 
   \qquad 
    \sigma\in \mathcal{R}_\alpha \Longrightarrow \dim_*\sigma \ge \alpha.
$$

Any finite Borel measure $\sigma \in \mathcal{M}^+(\mathbb{T})$ has the following Kahane decomposition.

 \begin{thm}[\cite{Kahane1988}] \label{thm:K-decomp} Every measure in  $\mathcal{M}^+(\mathbb{T})$ is uniquely decomposed into $\sigma_r + \sigma_s$
 with $\sigma_r \in \mathcal{R}_\alpha$ and $\sigma_s \in \mathcal{S}_\alpha$.
 \end{thm}  
 
 Actually the singular measure $\sigma_s$ is supported by the singular part of $\sigma$ defined by
 $$
   S_\sigma:= \left\{x: U^\sigma_\alpha(x) = \infty\right\}.
 $$
 The key step for proving  Theorem \ref{thm:K-decomp} is to prove ${\rm cap}_\alpha(S_\sigma)=0$.
 The regular measure $\sigma_r$ can then be decomposed as the sum of restrictions of $\sigma$ on $\{x: i-1\le U^{\sigma}_\alpha(x)<i\}$
 for $i\ge 1$. 
 \medskip

 It is clear that $\mathcal{S}_\alpha$ increases with $\alpha$, while $\mathcal{R}_\alpha$ decreases. Let 
 $$
       \mathcal{S}_{\alpha+} = \bigcap_{\beta >\alpha} \mathcal{S}_\beta, \qquad 
       \mathcal{R}_{\alpha -} = \bigcap_{\beta <\alpha} \mathcal{R}_\beta. 
 $$
 We can similarly define $\mathcal{S}_{\alpha-}$ and  $\mathcal{R}_{\alpha+}$:
 $$
       \mathcal{S}_{\alpha-} = \bigcup_{\beta <\alpha} \mathcal{S}_\beta, \qquad 
       \mathcal{R}_{\alpha +} = \bigcup_{\beta >\alpha} \mathcal{R}_\beta. 
 $$
 
 The following proposition, which is obvious,  shows how the potential theory can be used to  compute  the dimension of a measure. 
 
   \begin{prop}
  We have 
 $
      \sigma \in \mathcal{S}_{\alpha+} \cap  \mathcal{R}_{\alpha-} \Longrightarrow \dim \sigma =\alpha.
 $
 \end{prop}
 
 We also introduce the class $\mathcal{R}_\alpha^*$, which well describes the image of our operator $\mathbb{E}Q_\alpha$, as we shall prove. 
 A finite measure $\sigma \in \mathcal{R}_\alpha^*$ is defined by the following property:
 for any $\epsilon >0$, there exist a number $\beta>\alpha$ and  a compact set $K$ such that
 \begin{equation}\label{eq:R*}
      \sigma(K^c)<\epsilon;  \quad \forall x \in K, \ \ \underline{D}(\sigma1_K, x) \ge \beta.
 \end{equation}
 where $\sigma1_K$ means the restriction  to $K$ of $\sigma$. This class is not far from $\mathcal{R}_\alpha$.
 
 \begin{prop}
 We have the relations
 $
      \mathcal{R}_{\alpha+} \subset \mathcal{R}_\alpha^*\subset \mathcal{R}_\alpha.
 $
 \end{prop}
 
 \begin{proof} Assume $\sigma \in \mathcal{R}_{\alpha+}$. Then $\sigma\in \mathcal{R}_{\beta} $ for some $\beta>\alpha$, so 
 we can write $\sigma= \sum_{i=1}^\infty \sigma_i$ with $I^{\sigma_i}_\alpha <\infty$. Rewrite $\sigma=\sigma'+\sigma''$
 with $\sigma'=\sum_{i=1}^N \sigma_i$ and $\sigma''=\sum_{i=N+1}^\infty \sigma_i$.
 For any $\epsilon>0$, take $N$ large enough
 such that $\sigma''(\mathbb{T})<\frac{\epsilon}{2}$. Since $I^{\sigma'}_\alpha<\infty$, by Fubini theorem we have
 $U^{\sigma'}_\alpha(x)<\infty$ $\sigma'$-a.e. By the first estimate in  (\ref{eq:U-Lip}), we get 
 $\underline{D}(\sigma', x) \ge \beta$ $\sigma'$-a.e.  Recall that $\sigma'$ is a restriction of $\sigma$ on some set
 $B$. We can find a compact set $K$ in $B$ such that $\sigma(B\setminus K)<\frac{\epsilon}{2}$  so that 
 $\sigma(K^c)<\epsilon$, and $\underline{D}(\sigma1_K, x) \ge \beta$ for all $x\in K$. Thus we have proved that $\sigma \in \mathcal{R}_\alpha^*$.
 
 Now suppose that  $\sigma \in \mathcal{R}_\alpha^*$. Take $\epsilon_j =\epsilon 2^{-j}$, there exist a sequence  $\beta_j>\alpha$
 and compact set $K_j$ such that 
 $$
    \sigma(K_j^c)<\epsilon_j ;  \quad \forall x \in K_j, \ \ \underline{D}(\sigma1_{K_j}, x) \ge \beta_j >\alpha.
 $$
 Notice that $K= \bigcup K_j$ is large i.e. $\sigma(K^c)<\epsilon$. On the other hand, by using the second estimate in (\ref{eq:U-Lip}), 
we can prove that  $\sigma1_{K_j}\in \mathcal{R}_\alpha$ so that $\sigma 1_K \in \mathcal{R}_\alpha$. By letting $\epsilon \to 0$ , 
we can get other components of $\sigma$ which are of finite $\alpha$-energy.
 \end{proof}

    \section{Full action and $L^2$-theory }\label{sect:L2}
    Recall that our $\mathbb{T}$-martingale is defined by (\ref{eq:TM}) and the associated correlation function
    is 
    $
       H(t) = \frac{1}{2}\sum_{n=1}^\infty \alpha_n^2 \cos  nt. 
    $
    In this section we study the full action of the operator $Q$ defined by this $\mathbb{T}$-martingale with the second moment method.
    We shall come back to the investigation of  full action in the section \S \ref{sec:action} (cf. Theorem \ref{thm:ker}).
    The decomposition and composition discussed below are also preparing for this investigation.
    
    \subsection{$L^2$-convergence}
    For a large class of $\mathbb{T}$-martingales, there is a necessary and sufficient condition in terms of energy integral for the $L^2$-convergence 
    of the total mass martingale $\int Q_n(t)d\sigma$, which implies that $Q$ fully acts on $\sigma$.

       \begin{prop} \label{prop:L2}  Suppose  that $\alpha_n^2\downarrow 0$ and  there exists a constant $C$ such that
       \begin{equation}\label{eq:bd}
       \forall N\ge 1, \forall x\in \mathbb{T}, \quad 
      \frac{1}{2} \sum_{n=1}^N \alpha_n^2 \cos  nx\le H(x) +C.
    \end{equation}
    Let $\sigma$ be a probability measure on $\mathbb{T}$. 
       Then the martingale $\int Q_n(t) d\sigma(t)$ converges in $L^2$-norm if and only if the following energy is finite:
       $$
              \int\int e^{H(t-s)}d\sigma(t)d\sigma(s) <\infty.
       $$
   In this case, $\sigma$ is $Q$-regular.
     \end{prop}
   \begin{proof}
  The martingale $\int Q_n d\sigma$ converges in $L^2$-norm if and only if $\mathbb{E} \left(\int Q_nd\sigma \right)^2=O(1)$. But 
  $$
      \mathbb{E} \left(\int Q_nd\sigma \right)^2 =
      \int \int   \mathbb{E} [Q_n(t) Q_n(s)]d\sigma(t) d\sigma(s).
  $$
  By Lemma \ref{lem:EQQ}, we have
  $$
   \mathbb{E} \left(\int Q_nd\sigma \right)^2 \asymp 
       \int \int  \exp \left( \frac{1}{2}\sum_{k=1}^n \alpha_k^2 \cos  k (t-s) \right) d\sigma(t)d\sigma(s).
  $$
  Since the series defining $H$ convergence everywhere, we can conclude by using the hypothesis 
  (\ref{eq:bd})
   and the dominated convergence theorem.  
    \end{proof}
    
    We should point out that there is a difference between our trigonometric chaos and  other well studied chaos like the
     gaussian chaos \cite{Kahane1985}, the L\'evy chaos \cite{Fan1997}, the random covering 
    chaos \cite{Fan1989b}. The difference is that for these well studied chaos, 
    $\mathbb{E} Q_n(t)Q_n(s)$ increases with $n$, but it is not the case for the trigonometric chaos.  That is why we need the technical condition 
    (\ref{eq:bd}).

    \subsection{Energy of $Q\sigma$}
    When the random measure  $Q\sigma$ is non null, we would like to know if its $\alpha$-energy is finite.  
    We have the following result. Recall that $\Phi(t) = e^{H(t)}$.
       
  \begin{prop} \label{prop:chaos}Suppose that the kernel $\Phi$ satisfies the boundedness condition (\ref{eq:bd}) and suppose   
   $I_{\Phi \Phi_\alpha}^\sigma<\infty$. Then 
   we have $I^{Q\sigma}_{\alpha}<\infty$ a.s. Moreover 
  we have
  $$
      \mathbb{E} I^{Q\sigma}_{\alpha}  \asymp \int\int \Phi(t-s)\Phi_\alpha(t-s)d\sigma(t) d\sigma(s) <\infty.
  $$
  \end{prop}
  
\begin{proof}  
  Take an arbitrary number $L>0$. Denote
  $$
  \Phi_\alpha^{(L)}(t) =  L \wedge \Phi_\alpha(t),
  $$
    which is continuous and increases to  $\Phi_\alpha(t)$ as $L\uparrow \infty$. By the definition of $Q\sigma$ as weak limit and the continuity of 
    $\Phi_\alpha^{(L)}$, we have 
  $$
    \int \int \Phi_\alpha^{(L)}(t-s) dQ\sigma(t) dQ\sigma(s)   =  \lim_{n\to \infty} \int \int \Phi_\alpha^{(L)}(t-s) Q_n(t) Q_n(s) d\sigma(t) d\sigma(s).
  $$  
  By Fatou lemma, we get 
  $$
    \mathbb{E}\int \int \Phi_\alpha^{(L)}(t-s) dQ\sigma(t) dQ\sigma(s)  \le  \liminf_{n\to \infty} \int \int  \Phi_\alpha^{(L)}(t-s) \mathbb{E}Q_n(t) Q_n(s) d\sigma(t) d\sigma(s).
  $$  
  Replace $ \Phi_\alpha^{(L)}$ on the right hand side by $\Phi_\alpha$. Then apply the monotone convergence theorem on the left hand side, we get
  $$
    \mathbb{E} I^{Q\sigma}_\alpha \le  \liminf_{n\to \infty} \int \int \Phi_\alpha(t-s) \mathbb{E}Q_n(t) Q_n(s) d\sigma(t) d\sigma(s).
  $$  
  The estimate on $\mathbb{E}Q_n(t) Q_n(s)$ in Lemma \ref{lem:EQQ} and  the boundedness condition (\ref{eq:bd}) allow us to apply the dominated convergence theorem  to get    
$$
    \mathbb{E} I^{Q\sigma}_\alpha 
     \ll  \int \int   \Phi_\alpha(t-s)  \Phi(t-s)  d\sigma(t) d\sigma(s).
  $$ 
  
  To prove the inverse inequality, we start with
  $$
    \int \int \Phi_\alpha(t-s) dQ\sigma(t) dQ\sigma(s)   \ge   \liminf_{n\to \infty} \int \int \Phi_\alpha^{(L)}(t-s) Q_n(t) Q_n(s) d\sigma(t) d\sigma(s).
  $$
  The double integral on the right converges almost sure. It is bounded by $L (\int Q_n(t) d\sigma(t))^2$ which is uniformly integrable by Proposition \ref{prop:L2}.
  Therefore it converges in $L^1(\Omega)$, so that  
   $$
    \mathbb{E} I^{Q\sigma}_\alpha  \ge   \lim_{n\to \infty} \int \int  \Phi_\alpha^{(L)}(t-s) \mathbb{E} Q_n(t) Q_n(s) d\sigma(t) d\sigma(s).
  $$
  Again, by Lemma \ref{lem:EQQ} and the boundedness condition (\ref{eq:bd}),  we get 
  $$
   \mathbb{E} I^{Q\sigma}_\alpha  \gg  \int \int  \Phi_\alpha^{(L)}(t-s) \Phi(t-s) d\sigma(t) d\sigma(s).
  $$
  We finish the proof by letting $L\uparrow \infty$.
\end{proof}      

\subsection{Decomposition and composition of $Q_\alpha$} \label{subsect:decomp}
The following decomposition of $Q_\alpha$ will be useful when we study the regularity of $Q_\alpha\sigma$. Here the decomposition is in the sense of 
Theorem \ref{decomposition-principle} (a).

Define two sequences of weights as follows
$$
   P_n'(t) :=  P_{2n-1}(t), \qquad P_n''(t):= P_{2n}(t) 
$$
where $\{P_n(t)\}$ are the weights defined by (\ref{eq:W1}), with $\alpha_n = \frac{\alpha}{\sqrt{n}}$. 
The corresponding multiplicative chaotic operators of these two sequences of weights  will be denoted $Q'_{\alpha/\sqrt{2}}$ and $Q''_{\alpha/\sqrt{2}}$.
The reason for this parametrization $\alpha/\sqrt{2}$ is that the kernel associated to both $Q'_{\alpha/\sqrt{2}}$ and $Q''_{\alpha/\sqrt{2}}$
are the same as that of $Q_{\alpha/\sqrt{2}}$. 

Recall that  the correlation function of $Q_\alpha$ is equal to 
$$
         H_\alpha(t)= - \frac{\alpha^2}{2} \log \left|\sin \frac{t}{2}\right| +O(1).
   $$   
Notice that the correlation function of $Q''_\alpha$ is equal to 
\begin{equation}\label{eq:H1}
   H^{(2)} _\alpha(t) =\frac{\alpha^2}{2}\sum_{n=1}^\infty \frac{\cos 2n t}{2n} =  -\frac{\alpha^2}{4} \log \left|\sin t \right| - O(1),
\end{equation}
  and that of $Q''_\alpha$ is equal to 
  \begin{equation}\label{eq:H2}
   H^{(1)} _\alpha(t) =
   H_\alpha(t) - H^{(2)}_\alpha(t)  = -\frac{\alpha^2}{4} \log \left|\frac{\sin t/2}{2\cos t/2}\right| + O(1).
\end{equation}
Notice that $H^{(1)}$ and  $H^{(2)}$ have the same size as $H_{\alpha/\sqrt{2}}$ at $t=0$. After exponentiation, all three define the same 
Riesz kernel up to a multiplicative constant. 

\begin{figure}[htb]
 	\centering
 	\includegraphics[width=0.45\linewidth]{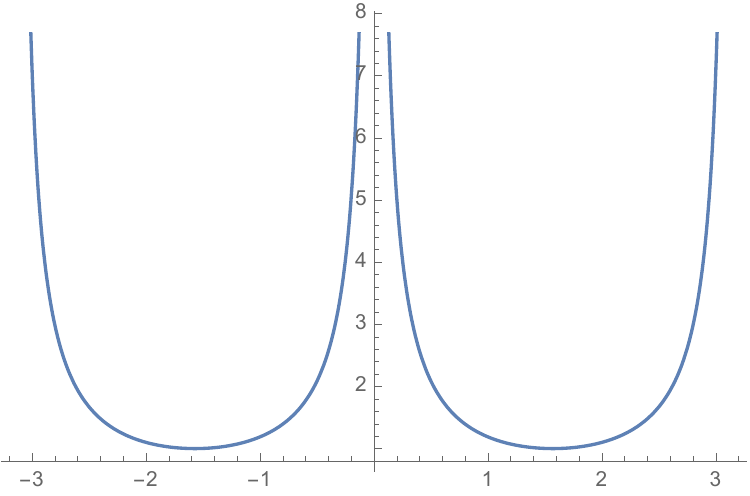}
 	\includegraphics[width=0.45\linewidth]{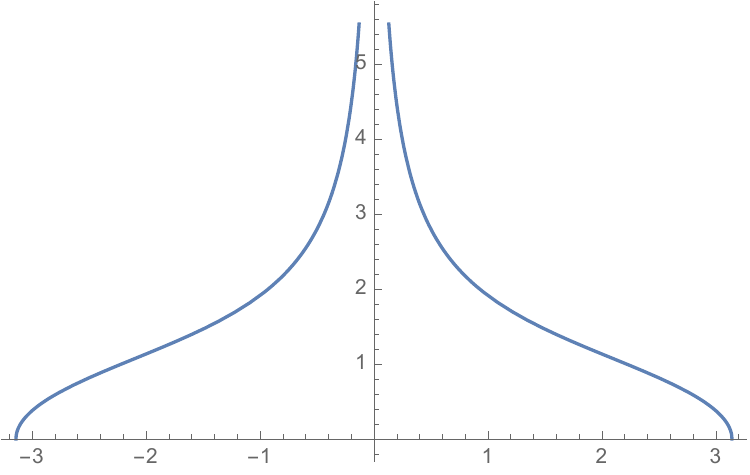}
 	\caption{The graphs of the kernels $e^{H_\alpha^{(1)}(t)}$ and $e^{H_\alpha^{(2)}(t)}$.}
 	\label{fig:1}
 \end{figure}

We are cheating on one point, but with no harm. 
   The kernel $e^{H_\alpha^{(2)}}$ has two singularities 
   $0$ and $\pi$ in the interval $[0, 2\pi)$. 
   We should consider that the operator $Q''_{\alpha/\sqrt{2}}$ is defined on $\mathbb{R}/\pi \mathbb{Z}$. 
   Actually the weights are $\pi$-periodic, not only $2\pi$-periodic. 
   The kernel $e^{H_\alpha^{(1)}}$ has only one singularity
   at $0$. 
   The figures of the kernels $e^{H_\alpha^{(1)}(t)}$ and $e^{H_\alpha^{(2)}(t)}$ are shown in Figure \ref{fig:1}. 
   
The decomposition principle stated in Theorem \ref{decomposition-principle} (a) gives rise to the
decomposition $Q_\alpha= Q''_{\alpha/\sqrt{2}} Q'_{\alpha/\sqrt{2}}$. We can continue in the same way  to decompose $Q'_{\alpha/\sqrt{2}}=
Q''_{\alpha/\sqrt{2}^2} Q'_{\alpha/\sqrt{2}^2}$, so that
$$
    Q_\alpha =Q''_{\alpha/\sqrt{2}} Q''_{\alpha/\sqrt{2}^2} Q'_{\alpha/\sqrt{2}^2}.
$$
The procedure can continue with $Q'_{\alpha/\sqrt{2}^2}$ and so on

Now let us talk about compositions. Let $Q_\alpha$ be the operator by the exponentiation of
$\frac{\alpha}{\sqrt{n}} \cos (nt+\omega_n')$ and   let $Q_\beta$ be the operator by the exponentiation of
$\frac{\beta}{\sqrt{n}}  \cos (nt+\omega_n'')$. We suppose that all $\omega_n'$ and $\omega_n''$ (for all $n\ge 1$)
are independent.  By the composition of $Q_\alpha$ and $Q_\beta$ we mean the operator by the exponentiation of 
all $\frac{\alpha}{\sqrt{n}}  \cos (nt+\omega_n')$ and $\frac{\beta}{\sqrt{n}}  \cos (nt+\omega_n'')$. In other word, it is the operator defined by the weights
$$
    P_{2n}(t) = \frac{e^{\frac{\alpha}{\sqrt{n}}  \cos (nt+\omega_n')}}{I_0(\frac{\alpha}{\sqrt{n}} )}, \quad P_{2n-1}(t) = \frac{e^{\frac{\beta}{\sqrt{n}}  \cos (nt+\omega_n'')}}{I_0(\frac{\beta}{\sqrt{n}}) }.
$$
It is clear that the correlation function of the composition operator is equal to
$$
          - \frac{\alpha^2+\beta^2}{2} \log \left|\sin \frac{t}{2}\right| +O(1).
   $$   
   This composition operator has the same potential properties as $Q_\gamma$ where $\gamma$ verifies
   $$
        \gamma^2 = \alpha^2 + \beta^2.
   $$

\subsection{Moments  with respect to Peyri\`ere measure}

  Let $\sigma$ be a probability measure on $\mathbb{T}$. If $Q$ acts fully on  $\sigma$, we define the Peyri\`ere  measure $\mathcal{Q}$ by the following relation
  $$
      \mathbb{E}_\mathcal{Q} \phi(t, \omega) = \mathbb{E} \int \phi(t, \omega) dQ\sigma(t).
  $$
 It can be proved that 
  $\{nt+\omega_n \mod 2\pi\}$ are $\mathcal{Q}$-independent random variables. This is a little more than what is stated in Theorem \ref{peyriere-thm}, but the proof is the same.
 The $\mathcal{Q}$-moment of $\varphi(nt +\omega_n)$ is easily computed, as stated in the following proposition. 
  

  \begin{prop} \label{prop:moments} Let $Q$ be the chaotic operator defined by $\{\alpha_n\}$ with $\alpha_n \to 0$. Suppose that $\sigma$ is a $Q$-regular probability measure.
  Let $\mathcal{Q}$ be the Peyri\`ere measure associated to $Q \sigma$. Then for any bounded or non-negative function $\varphi$ we have
  \begin{equation}\label{eq:ExpectQa}
       \mathbb{E}_\mathcal{Q} \varphi(nt + \omega_n) = \frac{1}{I_0(\alpha_n)} \int_0^{2\pi} \varphi(x)  e^{\alpha_n \sin  x} \frac{dx}{2\pi}.
  \end{equation}
  In particular,
   \begin{eqnarray*}
     \mathbb{E}_\mathcal{Q} \cos  (nt +\omega_n)
    & = &   \frac{I'_0(\alpha_n)}{I_0(\alpha_n)}
   = \frac{1}{2} \alpha_n+ O(\alpha_n^3)\\
     \mathbb{E}_\mathcal{Q} \cos^2  (nt +\omega_n) 
    & = &  \frac{I''_0(\alpha_n)}{I_0(\alpha_n)}
     = \frac{1}{2} + O(\alpha_n^2)
   \end{eqnarray*}
   \end{prop}
   
   \begin{proof} The formula (\ref{eq:ExpectQa}) is a special case of (\ref{eq:ExpectQ}).
   From  (\ref{eq:ExpectQa})  and Lemma \ref{lem:Bessel}, we get immediately the two other estimates.
   \end{proof}

  If we consider  $Q_\alpha$ defined by $\alpha_n =\frac{\alpha}{\sqrt{n}}$, we have
     \begin{eqnarray*}
     \mathbb{E}_\mathcal{Q} \cos  (nt +\omega_n)
    & = & 
    \frac{\alpha}{2}\cdot \frac{1}{\sqrt{n}} + O(n^{-3/2})\\
     \mathbb{E}_\mathcal{Q} \cos^2  (nt +\omega_n) 
    & = & 
      \frac{1}{2} + O(n^{-1})
   \end{eqnarray*}

   Consequently, we have 
   \begin{prop}  \label{prop:LogP} Make the same assumption as in Proposition \ref{prop:moments}.
   Assume further that $\sum \alpha_n^2 =\infty$. 
   Let $P_n(t)$ be the weights defining $Q$.
   Then for any $\delta>0$, almost surely for $Q\sigma$-almost every $t$
  we have 
  \begin{equation}\label{eq:LLN}
         \sum_{n=1}^N\log P_n(t) = \frac{1}{4} \sum_{n=1}^N \alpha^2_n + o\left( \left(\sum_{n=1}^N \alpha_n^2 \right)^{1/2+\delta} \right).
         \end{equation}
  \end{prop} 
  
  \begin{proof}  Using $\log I_0(\alpha_n) = \frac{1}{4}\alpha_n^2 + O(\alpha_n^4)$ and Proposition \ref{prop:moments}, we get
  $$
      \mathbb{E}_\mathcal{Q} \log P_n(t)= \alpha_n \mathbb{E}_\mathcal{Q} \cos (n t+\omega_n ) - \log I_0(\alpha_n) 
     = \frac{1}{4}\alpha^2_n + O(\alpha_n^4),
  $$
  and
  \begin{eqnarray*}
     &&\mathbb{E}_\mathcal{Q} (\log P_n(t))^2 \\  & = & \alpha_n^2 \mathbb{E}_\mathcal{Q} \cos^2 (n t+\omega_n ) - 2\alpha_n \log I_0(\alpha_n) 
      \mathbb{E}_\mathcal{Q} \cos (n t+\omega_n )  + \log^2 I_0(\alpha_n)\\
     &=& \frac{1}{4}\alpha^2_n + O(\alpha_n^4).
  \end{eqnarray*}
  So, we get the variance ${\rm Var}_{\mathcal{Q}}(\log P_n(t) )=  \frac{1}{4}\alpha^2_n + O(\alpha_n^4)$. 
  Let $Y_n = \log P_n(t)- \mathbb{E}_\mathcal{Q} \log P_n(t)$. Then the series 
  $$
       \sum_{n=1}^\infty \frac{Y_n}{[{\rm Var}_{\mathcal{Q}} (Y_1)+ \cdots + {\rm Var}_{\mathcal{Q}} (Y_n)]^{1/2+\delta}}
  $$
  converges $\mathcal{Q}$-a.e. Indeed, since $\{P_n\}$ are $\mathcal{Q}$-independent, the partial sums of the series is a $L^2$-bounded martingale
  and Doob's convergence theorem applies.  For the boundedness we use the fact that $\sum \frac{a_n}{(a_1+\cdots +a_n)^c}<\infty$
  for any positive numbers $a_n$ and $c>1$. 
  We conclude with the help of Kronecker lemma. 
  \end{proof}
  
   If we consider  $Q_\alpha$ defined by $\alpha_n =\frac{\alpha}{\sqrt{n}}$, we have
   \begin{equation}\label{eq:LLN2}
       \lim_{N\to \infty}   \frac{1}{\log N}  \sum_{n=1}^N\log P_n(t) =  \frac{\alpha^2}{4}.
         \end{equation}

  Intuitively, if $I$ is an interval  containing $t$ of length $N^{-1}$, we would have 
  $$
   Q_\alpha\sigma(I) \asymp  P_1(t)P_2(t)\cdots P_N(t) \sigma (I)
  $$
  so that 
  $$ 
     \lim_{|I|\to 0} \frac{\log  Q_\alpha\sigma(I)  }{\log  |I| } =   -   \lim_{N\to \infty}   \frac{1}{\log N}  \sum_{n=1}^N\log P_n(t) +   \lim_{|I|\to 0} \frac{\log  \sigma(I)  }{\log  |I| }.
  $$
   Therefore $\frac{\alpha^2}{4}$ is the difference of dimensions of the measure $\sigma$ and its image $Q_\alpha \sigma$:
   $$
        \dim Q_\alpha \sigma = \dim \sigma - \frac{\alpha^2}{4}.
   $$ 
   This will be rigorously proved (see Theorem \ref{thm:dim}).

\subsection{Mutual singularity and continuity}

Let $\{\alpha_n\}$ and $\{\alpha_n'\}$ be two sequences, which defines two operators $Q$ and $Q'$ using the same random variables
$\{\omega_n\}$.  There is a simple criterion for the mutual singularity and continuity of $Q\sigma$ and $Q'\sigma$.

\begin{prop} \label{prop:mutualsing}  Suppose $\sum |\alpha_n|^4 <\infty$ and $\sum |\alpha_n'|^4<\infty$. Suppose that  $\sigma$ is a $Q$-regular measure. We have
\begin{itemize}
  \item [(a)] $\sum_{n=1}^\infty |\alpha_n -\alpha_n'|^2<\infty \Rightarrow
      Q' \sigma \ll Q\sigma$ and $\sigma$ is $Q'$-regular.
  \item [(b)]  $\sum_{n=1}^\infty |\alpha_n -\alpha_n'|^2=\infty \Rightarrow
      Q' \sigma \perp Q\sigma$.
\end{itemize}
\end{prop}

\begin{proof}  This is a special case of Theorem \ref{ThmFan1991}. Let us compute
  $$
   \mathbb{E}\sqrt{P_n(t)P_n'(t)} = \frac{\mathbb{E} e^{ \frac{\alpha_n +\alpha_n'}{2} \cos (nt + \omega)}}{\sqrt{I_0(\alpha_n)I_0(\alpha_n')}}
   =   \frac{I_0(\frac{\alpha_n +\alpha_n'}{2})}{\sqrt{I_0(\alpha_n)I_0(\alpha_n')}}.
   $$
   Observe that 
   $$
     I_0(\frac{\alpha_n +\alpha_n'}{2}) = 1 + \frac{(\alpha_n +\alpha_n')^2}{16} + O(\alpha_n^4 + \alpha_n'^4),
   $$
   \begin{eqnarray*}
    I_0(\alpha_n)I_0(\alpha_n')
    & = &\left(1+\frac{\alpha^2_n}{4} + O(\alpha_n^4)\right) \left(1+\frac{{\alpha'}^2_n}{4} + O({\alpha'}_n^4)\right)\\
    & = & 1 +\frac{\alpha_n^2 +{\alpha'}_n^2}{4} + O(\alpha_n^4 + \alpha_n'^4).
   \end{eqnarray*}
   It follows that 
   \begin{eqnarray*}
         \mathbb{E}\sqrt{P_n(t)P_n'(t)} 
         & = & \exp\left( \frac{(\alpha_n +\alpha_n')^2}{16} -  \frac{\alpha_n^2 +{\alpha'}_n^2}{8} +   O(\alpha_n^4 + \alpha_n'^4) \right)\\
        &= & \exp\left( - \frac{(\alpha_n - \alpha_n')^2}{16}  +   O(\alpha_n^4 + \alpha_n'^4) \right).
   \end{eqnarray*}
   Finally we conclude by applying  Theorem \ref{ThmFan1991}.
\end{proof}

\section{Degeneracy }  
In this section we shall see that the chaotic operator $Q$ can be completely degenerate in the sense that 
$Q\sigma=0$ for all Borel probability measures $\sigma$.  
Recall that $H$, defined by (\ref{eq:H}), is the correlation function associated to $Q$. We denote its partial sums by  
$$
     H_n(t) = \frac{1}{2}\sum_{k=1}^n \alpha_n^2 \cos  kt.
$$

 \subsection{A general result}
 Let
 $$
     S_n(t) = \sum_{k=1}^n \alpha_k \cos  (kt+\omega_k).
$$
The following estimates will allow us to draw conditions for the complete degeneracy in different cases.  

\begin{prop} \label{prop:EstDeg} Let $0<h<1$ be a positive number,  $n\ge 1$ an integer, $I\subset \mathbb{T}$ an interval and 
$s\in I$ a fixed point. Let $(p, q)$ a conjugate pair such that $\frac{1}{p}+\frac{1}{q}=1$. We have 
\begin{equation}\label{eq:EST1}
\mathbb{E}\sup_{t\in I}Q_n(t)^h \le
\left(\mathbb{E} Q_n(s)^{hp}\right)^{1/p}  \left(\mathbb{E} e^{h|I| q \|S_n'\|_\infty} \right)^{1/q}; 
\end{equation}
\begin{equation}\label{eq:EST2}
 \left(\mathbb{E} Q_n(s)^{hp}\right)^{1/p}
 \ll  \left[ e^{ - \frac{1}{4}\cdot \frac{h(1 - hp)}{1-h} \sum_{k=1}^n \alpha_k^2}  \right]^{1-h};
\end{equation}
 \begin{equation}\label{eq:EST3}
 \left( \mathbb{E} e^{|I| hq \|S_n'\|_\infty} \right)^{1/q}
     \ll  n^{1/q} \exp\left( \frac{h^2}{4} q |I|^2 \sum_{k=1}^n k^2 \alpha_k^2 \right).
\end{equation}
\end{prop}

\begin{proof}
We write
$$
   Q_n(t) = Q_n(s) e^{ S_n(t)-S_n(s)} = Q_n(s) e^{(t-s) S_n'(r)} 
$$
where $r$ is a point between $s$ and $t$, and 
$$
S_n'(r) = - \sum_{k=1}^n k\alpha_k \sin  (kr +\omega_k).
$$
It follows that 
$$
   \sup_{t\in I}Q_n(t)^h \le  Q_n(s)^h e^{h|I| \|S_n'\|_\infty}. 
$$
Thus,  by H\"{o}lder's inequality, we get 
$$
\mathbb{E}\sup_{t\in I}Q_n(t)^h \le
\left(\mathbb{E} Q_n(s)^{hp}\right)^{1/p}  \left(\mathbb{E} e^{h|I|q \|S_n'\|_\infty} \right)^{1/q}, 
$$
which is (\ref{eq:EST1}).

The first expectation on the right hand side in the above inequality is easy to estimate by the independence and the local behavior of the Bessel function $I_0(\cdot)$ (cf. Lemma \ref{lem:Bessel}). Indeed, 
$$
 \mathbb{E} Q_n(s)^{hp}= \prod_{k=1}^n\frac{I_0(hp\alpha_k)}{I_0(\alpha_k)^{hp}}
 \ll e^{ - \frac{hp - h^2p^2}{4} \sum_{k=1}^n \alpha_k^2}.
$$
Here the assumption $\sum |\alpha_k|^4 <\infty$ is also used. Thus
\begin{equation*}
 \left(\mathbb{E} Q_n(s)^{hp}\right)^{1/p}
 \ll e^{ - \frac{h(1 - hp)}{4} \sum_{k=1}^n \alpha_k^2} = \left[ e^{ - \frac{1}{4}\cdot \frac{h(1 - hp)}{1-h} \sum_{k=1}^n \alpha_k^2}  \right]^{1-h},
\end{equation*}
which is (\ref{eq:EST2}).

Now let us prove the last desired  estimate. The arguments used below is inspired from \cite{Kahane1985}, p. 68-69).
Assume $\|S_n'\|_\infty = |S_n'(r)|$, where $r$ is a random point. By Bernstein's inequality, there exists an interval $J$ of length  $n^{-1}$ such that 
$|S'(\cdot)|\ge \frac{1}{2}|S'(r)|$. Thus
\begin{eqnarray*}
     e^{|I| hq \|S_n'\|_\infty} 
     &\le & \frac{1}{2 |J|} \int_J e^{|I| hq S_n'(x)|} dx
     \le   \frac{n}{2} \int_\mathbb{T} e^{|I| hq |S_n'(x)|} dx.
\end{eqnarray*}
Take the expectation to get
\begin{eqnarray*}
     \mathbb{E} e^{|I| hq \|S_n'\|_\infty} 
     &\le &   \frac{n}{2} \int_\mathbb{T}  \mathbb{E} e^{|I| hq |S_n'(x)|} dx.\\
     &\le &   \frac{n}{2} \int_\mathbb{T}\left( \mathbb{E} e^{|I| hq S_n'(x)} + \mathbb{E} e^{- |I| hq S_n'(x)} \right)dx.
 \end{eqnarray*}    
 The two last expectations don't depend on $x$, which can be easily and exactly computed. Thus
 \begin{eqnarray*}
  \mathbb{E} e^{|I| hq \|S_n'\|_\infty} 
  \le  n \prod_{k=1}^n I_0( |I| hq k\alpha_k)
     \ll n \exp\left( \frac{[ hq |I|]^2}{4} \sum_{k=1}^n k^2 \alpha_k^2 \right).
\end{eqnarray*}
Take $q$-th root of both sides to get (\ref{eq:EST3}) .
\end{proof}

\subsection{Degeneracy of $Q_\alpha$}

Let us apply Proposition \ref{prop:EstDeg} to treat the special case of $Q_\alpha$. 
\begin{thm} \label{thm:Degen}If $\dim^* \sigma <\frac{\alpha^2}{4}$, then $Q_\alpha \sigma =0$ a.s. 
In particular, 
the operator $Q_\alpha$ is completely degenerate when $|\alpha|>2$. 
\end{thm}

\begin{proof} Since $\alpha_n =\frac{\alpha}{\sqrt{n}}$, we have
$$
   \sum_{k=1}^n \alpha_k^2 \sim \alpha^2 \log n, \qquad   \sum_{k=1}^n k^2 \alpha_k^2 \sim \frac{\alpha^2}{2} n^2.
$$
Take a very small $\epsilon>0$. For a given interval $I$, choose $n$ such that $|I| \sim \frac{1}{n^{1+\epsilon}}$. 
Then choose $q= n^{\epsilon}$, so that $|I|^2 n^2 q = O(1)$.  Also notice that $n^{1/q} = n^{1/n^\epsilon}= O(1)$.
So, by Proposition \ref{prop:EstDeg}, we get
$$
 \mathbb{E} \sup_{t\in I} Q_n(t)^h
 \ll  \left(\frac{1}{n}\right)^{ \frac{\alpha^2}{4}\cdot \frac{h(1 - hp)}{1-h} (1-h)}
 \ll |I|^{ \frac{\alpha^2}{4}\cdot \frac{h(1 - hp)}{(1-h)(1+\epsilon)} (1-h)}.
$$
Since $\dim^* \sigma <\frac{\alpha^2}{4}$, 
 for $\epsilon >0$ small enough we have $\dim^* \sigma <(1-\epsilon) \frac{\alpha^2}{4}$. Take $1>h>1-\epsilon^2$.   Take $p$  close to $1$ (i.e. $|I|$ is sufficiently small) so that 
$$
   \frac{h(1 - hp)}{(1-h)(1+\epsilon)} \ge 1-\epsilon.
$$
We have proved that for any small interval $I$ (i.e. $|I|<\delta$ where $\delta$ depends on $\epsilon$ and $p$) there exists an integer $n$ such  that
$$
 \mathbb{E} \sup_{t\in I} Q_n(t)^h
 \ll |I|^{\frac{\alpha^2}{4}(1-\epsilon) (1-h)}.
$$
We conclude for the first assertion by Theorem \ref{thm:Kdegen}. The second assertion is a trivial consequence of the first one.
\end{proof}

If $\alpha_n = \frac{a_n}{\sqrt{n}}$ with $a_n \to +\infty$, the operator is completely degenerate.

\section{Image and kernel of  the projection $\mathbb{E}Q_\alpha$}
\label{sec:action}
Now we are going to show that 
if $\sigma$ is rather regular in the sense of Riesz potential theory, then it is $Q_\alpha$-regular. The following theorem gives an exact statement. 
This result and Theorem \ref{thm:Degen} on the degeneracy will give us a satisfactory description of the 
image and kernel of  the projection $\mathbb{E}Q_\alpha$. It is then possible to compute the dimension of 
$Q_\alpha \sigma$.

\subsection{$Q_\alpha$-regular measures}
\begin{thm}\label{thm:image} Suppose that $D >\frac{\alpha^2}{4}$. If $\sigma \in \mathcal{R}_D^*$, then $Q_\alpha \sigma \in \mathcal{R}_{D-\frac{\alpha^2}{4}}^*$.
 \end{thm}
\begin{proof}

We shall follow
  \cite{Kahane1985b}. The method was also used in \cite{Fan1997}.
  By the hypothesis on $\sigma$, for any $\epsilon>0$ there exist a number $\eta >D$  and a compact set $K$ such that $\sigma(K^c)<\epsilon$
  and 
      \[ I_\eta (\sigma 1_K) = \int_K \int_K \frac{d \sigma (t) d \sigma (s)}
         {| t - s |^\eta} < \infty. 
      \]
      
      {\em  Step I.  We first suppose that $D>\frac{\alpha^2}{2}$}. 
  By the  $L^2$-theory (Proposition \ref{prop:L2}),  $\mathbb{ E} Q_\alpha \sigma_K = \sigma_K$, where $\sigma_K$
  is the restriction $\sigma1_K$.
  The associated Peyri\`ere measure $\mathcal{Q}_K$ is well defined. 
  Denote by $R_n$
  the chaotic operator defined by the weights $P_m$ with $m \geq e^{n}$ and consider  the mass of $R_n\sigma_K$ 
  concentrated in the ball $B(t, e^{-n})$ centered at $t$ of radius $e^{-n}$:
      \[ \rho_n (t) = R_n \sigma_K (B(t, e^{-n})),
      \]
      which can be rewritten as 
       \[ \rho_n (t) = \int \chi_n (t, s) d R_n \sigma (s)
      \]
  where $\chi_n (t, s)$ is the characteristic function of the set of the points  $(t,s)\in K \times K$
  such that  $|t -s| \leq e^{-n}$.  We can then compute the $\mathcal{Q}_K$-expectation of $\rho_n$.  
  Indeed, by the definition of $\mathcal{Q}_K$ and Proposition \ref{prop:L2}, we have 
    \[ \mathbb{E}_{\mathcal{Q}_K} \rho_n  = \int \int \chi_n (t, s) \prod_{m=[e^{n}]}^\infty \mathbb{E}P_m (t) P_m (s)
         d \sigma_K (t) d \sigma_K (s).
    \]
  Notice that 
    \[ \sum_{m=1}^{[e^{n}]} \mathbb{E}P_m (t) P_m (s)  = \exp\left(\frac{\alpha^2}{2}\sum_{m=1}^{[e^n]} \frac{\cos m(t-s)}{m} + O(1) \right).
    \]
    If  $|t - s| \le e^{-n}$, then 
     \[ \sum_{m=1}^{[e^{n}]} \mathbb{E}P_m (t) P_m (s)  = \exp\left(\frac{\alpha^2}{2} n  + O(1) \right).
    \]
  Therefore
    \begin{eqnarray*}
       \mathbb{E}_{\mathcal{Q}_K} \rho_n  
        \ll
             \int \int e^{-\frac{\alpha^2}{2} n
             }\frac{\chi_n (t, s)}{\left|\sin\frac{t -s}{2}\right|^{\frac{\alpha^2}{2}}}
             d \sigma_K (t) d \sigma_K (s).
    \end{eqnarray*}
   If $ \chi_n (t, s) = 1$, we have   $e^n \le 1 / |t -s|$ so that 
    \[  \sum_{m=1}^\infty e^{(\eta - \frac{\alpha^2}{2}) n} \frac{\chi_n (t, s)}{\left|\sin \frac{t -s}{2}\right|^{\frac{\alpha^2}{2}}} \ll
             \frac{1}{|t-s|^\eta},
    \]
   It follows that 
    \[\mathbb{E}_{\mathcal{Q}_K} \left( \sum_{m=1}^\infty e^{\eta n} \rho_n \right) < \infty
    \]
  and consequently  $\rho_n = O(e^{-\eta n})$ $\mathcal{Q}_K$-a.s.  So, almost surely
    \[  \rho_n (t) = O(e^{-\eta n}) \,\, \,\, \quad Q_\alpha \sigma_K \! -\! a.e.
    \]
    Then for any $\epsilon > 0$, almost surely 
there exists a compact subset  $K_\epsilon^1\subset K$ such that
    \[ Q_\alpha \sigma (K\setminus K_\epsilon^1 ) < \epsilon
    \]
    \[ \limsup_{n \rightarrow \infty} \frac{\log \rho_n (t)}{n}
       \leq - \eta \,\,\, {\rm
      uniformly \ on } \, K_\epsilon^1.
    \]

  On the other hand,  from Proposition \ref{prop:LogP} we can deduce that  almost surely 
  for every $\epsilon > 0$, there exists a
  compact set $K_\epsilon^2$ contained in  $K$ such that 
     \[ Q_\alpha \sigma (K\setminus K_\epsilon^2 ) < \epsilon
    \]
    \[ \frac{\log Q_n (t)}{\log n}
       \rightarrow \frac{\alpha^2}{4} \,\,\, {\rm
       uniformly \, on} \, K_\epsilon^2.
    \]

  Now let
    \[ S_\epsilon = 1_{K_\epsilon^1 \bigcap K_\epsilon^2} Q \sigma.
    \]
  Then
     \[ Q \sigma (K \setminus (K_\epsilon^1 \bigcap K_\epsilon^2)) < 2 \epsilon,
     \]
and
  \begin{eqnarray*}
         \limsup_{n \rightarrow \infty} \frac{\log S_\epsilon
            (B (t, e^{-n}))}{n} 
        \leq 
        \limsup_{n \rightarrow \infty} \frac{\log Q_{[e^{dn}]}(t)}
        {n} + \limsup_{n \rightarrow \infty} \frac{\log \rho_n (t)}
        {n} 
       \leq 
         \frac{\alpha^2}{4}- \eta
  \end{eqnarray*}
  uniformly on $ K_\epsilon^1 \bigcap K_\epsilon^2 $. This implies
   $Q_\alpha \sigma \in M_{(D - \frac{\alpha^2}{4})^+}^+ ({\mathbb{T}})$. Let us summarize as follows
   \begin{equation}\label{eq:decomp-principle}
      D >\frac{\alpha^2}{2}, \ \ \sigma \in \mathcal{R}^*_D \Longrightarrow Q_\alpha\sigma_K \in \mathcal{R}^*_{D-\frac{\alpha^2}{4}},
   \end{equation}
   where $\sigma(K^c)$ can be made as small as we like.
   
   {\em Step II. General case.} We use the decomposition $$
   Q_\alpha = Q^{(m)}_{\frac{\alpha}{\sqrt{2}^m}} Q'_{\frac{\alpha}{\sqrt{2}^m}} \cdots Q'_{\frac{\alpha}{\sqrt{2}^2}} Q'_{\frac{\alpha}{\sqrt{2}}}$$
   which is  presented in 
   \S \ref{subsect:decomp}. First apply the principle (\ref{eq:decomp-principle})
 to $Q'_{\frac{\alpha}{\sqrt{2}}}$  and $\sigma$. We check that $D>\frac{\alpha^2}{4}=\frac{(\alpha/\sqrt{2})^2}{2}$
 and $\sigma  \in \mathcal{R}^*_D$, and we get 
 $$
 Q'_{\frac{\alpha}{\sqrt{2}}} \sigma \in  \mathcal{R}^*_{D-\frac{\alpha^2}{4\cdot 2}}.
 $$
  Now observe that $D-\frac{\alpha^2}{4\cdot 2}> \frac{\alpha^2}{8}=\frac{(\frac{\alpha}{\sqrt{2}^2})^2}{2}$. So, the 
   principle (\ref{eq:decomp-principle}) applies 
  $Q'_{\frac{\alpha}{\sqrt{2}^2}}$  and $Q'_{\frac{\alpha}{\sqrt{2}}} \sigma$. Thus we get 
  $$
   Q'_{\frac{\alpha}{\sqrt{2}^2}} Q'_{\frac{\alpha}{\sqrt{2}}} \sigma \in  \mathcal{R}^*_{D-\frac{\alpha^2}{4\cdot 2}-  \frac{\alpha^2}{4\cdot2^2}}.
  $$
  Inductively we get 
  $$
   Q'_{\frac{\alpha}{\sqrt{2}^m}}  \cdots  Q'_{\frac{\alpha}{\sqrt{2}^2}} Q'_{\frac{\alpha}{\sqrt{2}}} \sigma \in  \mathcal{R}^*_{D-\frac{\alpha^2}{4\cdot 2}-  \frac{\alpha^2}{4\cdot 2^2}
   - \cdots-\frac{\alpha^2}{ 4\cdot 2^m}}.
  $$
 Since $D-\frac{\alpha^2}{4}>0$. When $m$ is sufficiently large, we have 
 $$
     D-\frac{\alpha^2}{4\cdot 2}-  \frac{\alpha^2}{4\cdot 2^2}
   - \cdots-\frac{\alpha^2}{ 4\cdot 2^m} > \frac{(\alpha/\sqrt{2}^m)^2}{2}= \frac{\alpha^2}{2^{m+1}}.
 $$
 So, 
 we can apply once more the
   principle (\ref{eq:decomp-principle}), to  $Q^{(m)}_{\frac{\alpha}{\sqrt{2}^m}} $. Finally we get 
  $$
  Q^{(m)}_{\frac{\alpha}{\sqrt{2}^m}}  Q'_{\frac{\alpha}{\sqrt{2}^m}}  \cdots  Q'_{\frac{\alpha}{\sqrt{2}^2}} Q'_{\frac{\alpha}{\sqrt{2}}} \sigma \in  \mathcal{R}^*_{D-\frac{\alpha^2}{4\cdot 2}-  \frac{\alpha^2}{4\cdot 2^2}
   - \cdots-\frac{\alpha^2}{ 4\cdot 2^m} -\frac{\alpha^2}{ 4\cdot 2^m}} = \mathcal{R}^*_{D-\frac{\alpha^2}{4\cdot 2}}.
  $$ 
  
  \end{proof}
  
  \subsection{${\rm Ker}\, \mathbb{E}Q_\alpha$ and ${\rm Im}\, \mathbb{E}Q_\alpha$}
  
  We have the following description  for the kernel and the image   $\mathbb{E} Q_\alpha$.
  
  \begin{thm} \label{thm:ker} We have 
  $$
  \mathcal{S}_{\frac{\alpha^2}{4}-} \subset {\rm Ker} \mathbb{E} Q_\alpha \subset \mathcal{S}_{\frac{\alpha^2}{4}+};
  $$
   $$
  \mathcal{R}_{\frac{\alpha^2}{4}+} \subset {\rm Im} \ \mathbb{E} Q_\alpha \subset \mathcal{R}_{\frac{\alpha^2}{4}-}.
  $$
  \end{thm}
  \begin{proof} The first inclusion follows immediately from Theorem \ref{thm:Degen}. 
  Now we prove the second inclusion by contradiction. Suppose $\sigma \in {\rm Ker} \mathbb{E} Q_\alpha$ but 
  $\sigma \not\in \mathcal{S}_{\frac{\alpha^2}{4}+}$. Then $\sigma \in \mathcal{S}^c_\beta$ for some $\beta > \frac{\alpha^2}{4}$.
  That means $\sigma$ has its $\beta$-regular component $\sigma_r \not =0$ according to the Kahane decomposition.
  Hence, by Theorem \ref{thm:image}, 
  $$\dim_*Q_\alpha \sigma \ge \dim_*Q_\alpha \sigma_r \ge \beta- \frac{\alpha^2}{4}>0 \quad  a.s
  $$ 
 which contradicts   $\sigma \in {\rm Ker} \mathbb{E} Q_\alpha$.
 
 Assume $\sigma \in  \mathcal{R}_{\frac{\alpha^2}{4}+}$. That means $\sigma \in  \mathcal{R}_\beta$ for some $\beta>{\frac{\alpha^2}{4}}$.
 We claim that $\sigma\in {\rm Im} \ \mathbb{E} Q_\alpha$. Otherwise, there is a non null component of $\sigma$, say $\sigma1_K$ with $\sigma(K)>0$, which is killed by
 $\mathbb{E} Q_\alpha$. That is to say, $\sigma 1_K \in {\rm Ker} \mathbb{E} Q_\alpha$.  So, $\sigma 1_K \in \mathcal{S}_{\frac{\alpha^2}{4}+}$,
 by what we have just proved above. Hence  ${\rm Cap}_\beta (K) =0$. This contradicts $\sigma 1_K \in \mathcal{R}_\beta$.
 
 The last inclusion is also proved by contradiction. Assume $\sigma \in \mathbb{E} Q_\alpha$ but  $\sigma \not \in \mathcal{R}_{\frac{\alpha^2}{4}-}$.
 Then $\sigma  \not\in  \mathcal{R}_\beta$ for some $\beta <\frac{\alpha^2}{4}$. By the Kahane decomposition, $\sigma$ has a non null $\beta$-singular component
 which must be killed by $\mathbb{E} Q_\alpha$ because $\beta <\frac{\alpha^2}{4}$ (cf. Theorem \ref{thm:Degen}). This contradicts $\sigma \in \mathbb{E} Q_\alpha$.
 \end{proof}
  
  Recall that $\dim \sigma=D$ means $\dim_*\sigma =\dim^* \sigma=D$.  In this case, we say that
  $\sigma$ is unidimensional.  The dimension of $Q_\alpha \sigma$ is given by the following formula.

\begin{thm} \label{thm:dim} 
Suppose that  $\sigma$ is a unidimensional measure and 
$\dim \sigma >\frac{\alpha^2}{4}$. We have $$\dim Q_\alpha \sigma = \dim \sigma - \frac{\alpha^2}{4}.$$
\end{thm}
\begin{proof} Let $D=\dim \sigma$. For any $\epsilon >0$, the fact $\dim_* \sigma=D$ implies $\sigma \in \mathcal{R}^*_{D-\epsilon}$. Then by Theorem \ref{thm:image}, we get 
$\dim_* Q_\alpha \sigma \ge D-\epsilon -\frac{\alpha^2}{4}$ a.s. This proves $\dim_* Q_\alpha\sigma \ge D-\frac{\alpha^2}{4}$ a.s.

Take $\beta$ such that $\frac{\beta^2}{4} +\frac{\alpha^2}{4}  >D$. Let $\gamma$ be such that $\gamma^2 =\beta^2 +\alpha^2$. 
We consider an operator $Q_\beta$, which is independent of $Q_\alpha$. The product $Q_\beta Q_\alpha$ has the same law as 
$Q_\gamma$.
The fact $\dim^* \sigma=D$, together with $\frac{\beta^2}{4} +\frac{\beta^2}{4}  >D$,  implies $\sigma \in {\rm Ker}\mathbb{E} Q_\gamma$. 
Then, by Theorem \ref{decomposition-principle} c), we get $Q_\alpha\sigma \in {\rm Ker} \mathbb{E}Q_\beta$ a.s.   
Now we apply Theorem \ref{thm:ker} to get $\dim^* Q_\alpha \sigma \le \beta$ a.s. Optimizing $\beta$ leads to
$\dim^* Q_\alpha \sigma \le D -\frac{\alpha^2}{4}$ a.s. 
\end{proof}
     

\section{Trigonometric multiplicative chaos on $\mathbb{T}^d$}
Our theory can be extended to torus $\mathbb{T}^d$ of all dimensions, with few points to be checked in higher dimension case.
On  $\mathbb{T}^d$ we propose to study the random series
\begin{equation}\label{eq:RSd}
 \sum_{n\in \mathbb{Z}^d_+}\rho_n \cos  (n\cdot x +\omega_n)
\end{equation} 
where $\mathbb{Z}^d_+$ is the set $\mathbb{Z}^d\setminus \{0\}$ but with $n$ and $-n$ identified, and $\{\omega_n\}$ are  iid random variables
uniformly distributed on $\mathbb{T}$ (not on $\mathbb{T}^d$). It can be rewritten as
$$
     \frac{1}{2}\sum_{n \in \mathbb{Z}^d \setminus \{0\}} \rho_n e^{i (n\cdot x+ \omega_n)}
$$
with $\rho_n = \rho_{-n}$ and $\omega_n =\omega_{-n}$. When, $d=1$, $\mathbb{Z}^d_+$ is the set of natural numbers. The correlation function
of (\ref{eq:RSd}) is defined by 
$$
 H (t) = \frac{1}{2}\ \sum_{n\in \mathbb{Z}^d_+}\rho_n^2 \cos  (n\cdot t) \qquad (t \in \mathbb{T}^d). 
$$

Let $|x|_\infty= \max_{1\le i\le d}|x_j|$ for $x =(x_1, \cdots, x_d)\in \mathbb{R}^d$.
Let $|x|$ to denote the euclidean norm $\sqrt{\sum_{i=1}^n x_i^2}$.

\subsection{Definition of $\mathbb{T}^d$-martingales}

Let $\{\alpha_k\}_{k\in \mathbb{Z}^d_+}$ be a family of real numbers such that $\sum_{k \in \mathbb{Z}^d_+ } |\alpha_k|^4 <\infty$.
 Let $\{A_n\}$ be an increasing sequence of finite subsets of $\mathbb{Z}^d_+$ such that $A_n \uparrow \mathbb{Z}^d_+$ with $A_0=\emptyset$. 
We define the weights
 \begin{equation}\label{eq:weightd}
       P_n(t) = \frac{\exp \sum_{k\in  A_n\setminus A_{n-1}} \alpha_k \cos (k \cdot t +\omega_k)}{ \prod_{k\in A_n\setminus A_{n-1}} I_0(\alpha_k)}.
 \end{equation}
 Here is a choice for $\{A_n\}$.
For $n\ge 1$, we consider the hypercube $C_n = \{k \in \mathbb{Z}^d: |k|_\infty \le  n\}$ and choose $A_n$ to the positive part 
$ C_n^+ = C_n \cap \mathbb{Z}^d_+$. In this case $ C_n^+\setminus C_{n-1}^+$ is the boundary $ \partial C_n^+$.
We have
$
   \#C_n = (2n+1)^d-1
$.
It follows that 
$$
 \#C_n^+ = \frac{(2n+1)^d-1}{2}, \qquad \#\partial C_n^+ 
                     \sim d 2^{d-1} n^{d-1}.
 $$
 Another choice for  $A_n$ is the ball $\{k\in \mathbb{Z}^d_+: |k|\le n\}$. 
 
 The $\mathbb{T}^d$-martingale defined by the weight (\ref{eq:weightd})  produces
a chaotic operator $Q$.  The same computation as in the case of dimension $1$ leads to the basic relation
 \begin{equation}\label{eq:fundamentd}
    \mathbb{E} Q_n(t) Q_n(s) = \exp \left(\frac{1}{2}\sum_{k\in A_n} \alpha_k^2 \cos k\cdot (t-s) + O(1) \right). 
 \end{equation}
 
 For a given real number $\alpha$, we have the typical sequence of coefficients
 $$
         \alpha_k = \frac{\alpha}{|k|^{d/2}}.
 $$
As we shall see below (cf. Theorem \ref{thm:G}),  the corresponding correlation function $ H_\alpha(t) $ is smooth in $\mathbb{T}^d \setminus \{0\}$ and around $0$ it behaves as
$$
     H_\alpha(t) = \frac{\alpha^2}{2}\sum_{k\in \mathbb{Z}^d_+} \frac{\cos k\cdot t}{|k|^d} \sim \frac{\alpha^2}{2} \cdot \frac{\pi^{d/2}}{\Gamma(d/2)} \log \frac{1}{|t|} 
     \quad  (t \to 0). 
 $$
 Its exponentiation  $\Phi_\alpha = e^{H_\alpha(t)}$ is then a Riesz kernel
 \begin{equation}\label{eq:Rieszd}
   \Phi_\alpha(t) \asymp \frac{1}{|t|^{\tau(\alpha) \alpha^2/2}} \qquad {\rm with} \ \ \  \tau(\alpha) = \frac{\pi^{d/2}}{\Gamma(d/2)}.
 \end{equation}

\subsection{Correlation function $H_\alpha$}
We consider the Jacobi function $G$ defined by the trigonometric series 
\begin{equation}\label{eq:G}
G(x)=\sum_{k\in {\mathbb Z}^d, k\neq 0}|k|^{-d}\exp(i k\cdot x).
\end{equation}
For $m\ge 1$, we denote by $S_m$ its  partial sum:
  $$
  S_m(G)=\sum_{ |k|\leq m, k\neq 0}|k|^{-d}\exp(i k\cdot x).
  $$
  We set $s_d= 2\frac{\pi^{d/2}}{\Gamma(d/2)} $ which is the area of the unit sphere in $\mathbb{R}^d$.

\begin{thm}\label{thm:G} The Jacobi function $G$ has the following properties:\\
\indent {\rm (a)}\ 
The function $G(x)$ is ${\mathcal C}^{\infty}$ in ${\mathbb T}^d\setminus \{0\}.$ \\
\indent  {\rm (b)}\  There is a bounded function $E(x)$ on  $\mathbb{T}^d$ such that 
 $$\forall x\in \mathbb{T}^d, \ \ \ \ G(x)=s_d\log\frac{1}{|x|}+E(x).$$ 
 
\indent {\rm (c)}\ 
There exists a constant $C$ such that for any integer $m\ge 1$  we have 
$$ 
S_m(G)(x)\leq s_d \log\frac{1}{|x|} +C.
$$
\end{thm}

\begin{proof} The proof mimics the argument given by Titchmarsh and Riemann himself of the functional equation of the zeta function.   One considers
the Jacobi  function $\theta(u, x)$ defined by  
 $$\theta(u, x)=\sum_{n \in \mathbb{Z}^d}  \exp(-|n|^2u)\exp(i n \cdot x).$$
  Here $u>0$ and $x\in \mathbb{T}^d$. The notations used here slightly differ from the standard ones. 
   Let us consider $\theta_0=\theta-1$. Then the function  $\theta_0(u, x)$ has an exponential decay as  $u$ tends to $+\infty$ and this decay is uniform in $x$. 
   Assume $s>d$. We introduce the auxiliary Jacobi functions 
  $$
  G(s, x)=\sum_{n \in \mathbb{Z}^d, n\neq 0}|n|^{-s}\exp(i n\cdot x).
  $$ 
  where the series converges absolutely. Fubini theorem implies
$$
 \int_0^\infty \theta_0(u, x)u^{s/2}\frac{du}{u}= \Gamma(s/2)G(s, x).
$$
One splits  this  integral into $\int_0^1+\int_1^\infty$ which yields
 $$
\Gamma(s/2) G(s, x)=A(s, x)+B(s, x).$$ 
The function $B(s, \cdot)$ is obviously a $C^\infty$ function for every $s\geq 0$, especially $B(d, x)$.
To study $A(s, x)$, one rewrites
$\theta_0=\theta-1$ and uses the following functional equation for the Jacobi function 
$$
\theta(u, x)=\frac{\pi^{d/2}}{u^{d/2}}\sum_{k \in \mathbb{Z}^d} \exp\Big(-\frac{|x- 2\pi k|^2}{4u}\Big),
$$ 
which follows from the Poisson summation formula.
One can ignore  $1$ in $\theta_0=\theta-1$ since its contribution to  $A(s,x)$ yields a constant. 
When $x$ belongs to the ball centered at $0$ with radius $1$  all the terms together in the above sum but the term $k=0$ yield a  ${\mathcal C}^{\infty}$ contribution. Therefore one is led to 
$$J(s,x) :=\pi^{d/2} \int_0^1 \exp\left(-\frac{|x|^2}{4u}\right)u^{\frac{s-d}{2}}\frac{du}{u}.$$ 
We have thus proved the following lemma.

\begin{lem}
The difference  \begin{equation}\label{eq:R}
\Gamma(s/2)G(s, x)-J(s, x)=R(s, x)
\end{equation} is a ${\mathcal C}^{\infty}$ function for every  $s\geq 0$. 
\end{lem}

Therefore one can pass to the limit  $s\to d$ in (\ref{eq:R}) and we are led to the computation of $J(d, x)$.
 Make a change of variables 
  $t= \frac{|x|^2}{4} u$ we get
  $$
       J(d, x) = \pi^{d/2}\int_0^{\frac{4}{x^2}} e^{-\frac{1}{t}}
 \frac{dt}{t}  = 2 \pi^{d/2} \log \frac{1}{|x|} +O(1), \qquad x \to 0.
 $$
 We have thus proved (a) and (b).

(c) will follow from (a), (b) and the  following lemma. 

\begin{lem}\label{lem:partials}
If $\phi$ denotes a non negative compactly supported radial function of integral $1$ and if $\phi_m(x)=m^d\phi(mx)$ we have 
$$\|G*\phi_m-S_m(G)\|_\infty\leq C.$$ 
\end{lem}
The proof of this lemma is immediate. We have $$
 |{\widehat \phi}\big(\frac{y}{m}\big)-1|\leq C\frac{|y|^2}{m^2}\ \ {\rm if } \ \ |y|\leq m;
 \quad |{\widehat \phi}\big(\frac{y}{m}\big)|\leq C\frac{m^2}{|y|^2} \ \ {\rm  if}\ \  |y|\geq m.
 $$ This together with the decay of the Fourier coefficients of $G$ yields the result. 

Now prove (c). By (a) and (b), the function $G$ can be forgotten and replaced by $f(x)=s_d|\log|x||$, because the question arises   around $0$.  It suffices to check that $f*\phi_m\leq f+C$ which is an easy calculation. 
 \end{proof}
 
 It is easy to deduce that the partial sums over cubes $C_n$ are  also bounded by $s_d \log (1/|x|) +C$.

  \subsection{Riesz potential theory on $\mathbb{T}^d$ and $L^2$-theory}
  Let $\beta\in (0, d)$ be a fixed number. It is known that (cf. \cite{SW1971}, p. 256)
  $$
      \sum_{m\in \mathbb{Z}^d, m\not=0} \frac{1}{|m|^{d-\beta}} e^{ i m\cdot x} \sim \gamma_\beta \frac{1}{|x|^{\beta}} + b(x)
  $$
  where $\gamma_\beta = 2^\beta \pi^{2\beta- \frac{d}{2}}\frac{\Gamma(\frac{d-\beta}{2})}{\Gamma(\frac{\beta}{2})}$ and $b\in C^\infty(\mathbb{T}^d)$.
  This means that the function on the right, which is  Lebesgue integrable on $\mathbb{T}^d$,  admits the series on the left  as its Fourier series.
   The function $b$ is necessarily real valued. 
  Let $B=\max_x b(x)$. 
  We define the $\beta$-order Riesz kernel by
  $$
     R_\beta(x) = \gamma_\beta \frac{1}{|x|^{\beta}} + b(x) +B.
  $$
  It is a positive, lower semicontinuous  function having positive Fourier coefficients. We also have the following simple estimate
  $$
      R_\beta(x) \asymp \frac{1}{|x|^\beta}   \qquad (|x_1|\le \pi, \cdots, |x_d| \le \pi).
  $$
  
  As in the case $d=1$, the $\beta$-order energy $I^\sigma_\beta$ of a measure $\sigma$ on $\mathbb{T}^d$ can be defined. As shown in  \cite{KS1963}, we have the formula
  $$
         I^\sigma_\beta = B|\widehat{\sigma}(0)|^2 + \sum_{|m|>0} \frac{|\widehat{\sigma}(m)|^2}{|m|^{d-\beta}}.
  $$
  More generally we have the results on $\mathbb{T}^d$ as developed in \S \ref{sect:capacity} about  potential, capacity and dimension
  (cf. \cite{KS1963},\cite{Kahane1988}, \cite{Fan1994}).
  
  By Theorem \ref{thm:G}, we have
        $$
            H_\alpha(t) = \frac{\alpha}{2}\cdot \frac{\pi^{d/2}}{\Gamma(d/2)} \log \frac{1}{|t|}+ O(1) \quad ({\rm as}\ \  t \to 0).
        $$
         Theorem \ref{thm:G} allows us to develop 
         the
   $L^2$-theory for our $\mathbb{T}^d$-martingales, which was developed in \S \ref{sect:L2} for the case $d=1$.  We don't restate the results here,  but we are free to use them.

  \subsection{Degeneracy}
  
  \begin{prop}
  The result of Proposition \ref{prop:EstDeg} remain true when $d\ge 2$. It suffices to replace 
  $\sum_{k=1}^n \alpha_k^2$ and $\sum_{k=1}^n k^2\alpha_k^2$ respectively by 
  $$
      \sum_{k\in A_n}\alpha_k^2, \qquad    \sum_{k\in A_n} |k|^2\alpha_k^2.
  $$
  \end{prop}
  \begin{proof} The proof is the same as that of Proposition \ref{prop:EstDeg}.
  We only need to point out one point. Let 
  $$
       S_n(t) = \sum_{k\in A_n} \alpha_k \cos (k\cdot t +\omega).
  $$
  Now
  $$
     |S_n(t) - S_n(s)| = \left|\sum_{j=1}^d \frac{\partial S_n}{\partial t_j}(r) (t_j-s_j)\right|
     \le |t-s|_\infty  \sum_{j=1}^d \left|\frac{\partial S_n}{\partial t_j}(r) \right|,
  $$
  where
  $$
   \frac{\partial S_n}{\partial t_j}(r)  =  - \sum_{k\in A_n} k_j \alpha_k \sin (k\cdot t +\omega)
  $$
  Let $I$ be a cube and let $|I|$ be its side length. We would have to estimate the following expectation
  $$
     \mathbb{E} \exp\left( h |I|q \sum_{j=1}^d\left\|\frac{\partial S_n}{\partial t_j}\right\|_\infty\right),
  $$
  which, by H\"{o}lder inequality, is bounded by 
  $$
      \left\{ \prod_{j=1}^d  \mathbb{E} \exp\left( d h |I|q \left\|\frac{\partial S_n}{\partial t_j}\right\|_\infty\right)\right\}^{1/d}
  $$
  Each of these expectation can be estimated as before. Finally we get 
  $$
   \mathbb{E} \exp\left( h |I|q \sum_{j=1}^d\left\|\frac{\partial S_n}{\partial t_j}\right\|_\infty\right)
   \ll n \exp\left( \frac{d [ hq |I|]^2}{4} (k_1^2+\cdots + k_d^2) \alpha_k^2 \right).
   $$
   \end{proof}
  \medskip
  
  Let us apply the result to the special case $\alpha_k = \frac{\alpha}{|k|^{d/2}}$. Let us take 
  $A_n =C_n^+$. We are led to estimate
  $$
       \sum_{k\in C_n^+} \alpha_k^2 =  \sum_{k\in C_n^+} \frac{\alpha^2}{|k|^d}, 
        \qquad \sum_{k\in C_n^+} |k|^2\alpha_k^2 =  \sum_{k\in C_n^+} \frac{\alpha^2}{|k|^{d-2}}.
  $$
  The sum over $C_n$ doubles the sum over $C_n^+$. The sum over $C_n$ and the sum over the ball $B(0,n)$ are comparable, because  
 the cube $C_n$ is contained in the ball centered at $0$ of radius $\sqrt{d} n$, and contains the ball of radius $n$.

  \begin{lem} Let $B(0, R)$ be the ball of radius $R$ centered at $0$. We have
  $$
  \sum_{1\le |k|\le R} \frac{1}{|k|^d} =  s_d \log R + O(1), 
  \qquad  \sum_{1\le |k|\le R} \frac{1}{|k|^{d-2}} =s_d R^2 + O(R).
  $$
  where $s_d= 2 \frac{\pi^{d/2}}{\Gamma(d/2) }$. Also we have
  $$
      \int_{n\le |x|\le \sqrt{d}n}\frac{dx}{|x|^d} = O(1), \qquad  \int_{n\le |x|\le \sqrt{d}n}\frac{dx}{|x|^{d-2}} \asymp n^2.
  $$
  \end{lem}
  \begin{proof} Assume $d\ge 2$. When $x$ is in the unit cube centered by $k$, we have $|k|^{-d}- |x|^{-d} = O(|k|^{-(d+1)})$. 
  Since  $\sum |k|^{-(d+1)}<\infty$, we have
  $$
       \sum_{1\le |k|\le R} \frac{1}{|k|^d} = \int_{1\le |x|\le R}\frac{dx}{|x|^d} + O(1).
  $$
 The integral is equal
  to $s_d \int_1^R \frac{dr}{r} = s_d \log R$.  The second equality can be proved in the same. We need the estimate
   $|k|^{-d+2}- |x|^{-d+2} = O(|k|^{-(d-1)})$ for $x$ in the unit cube centered by $k$.
   The estimates of two integrals are direct. 
  \end{proof}
  
  This lemma implies
    \begin{eqnarray*}
       \sum_{k\in C_n^+} \alpha_k^2 & =&  \sum_{k\in C_n^+} \frac{\alpha^2}{|k|^d} = \frac{\alpha^2}{2} s_d \log n + O(1),\\
        \qquad \sum_{k\in C_n^+} |k|^2\alpha_k^2 & =&   \sum_{k\in C_n^+} \frac{\alpha^2}{|k|^{d-2}} = \frac{\alpha^2}{2} s_d n^2 + O(n).
  \end{eqnarray*}
  
  Now we are ready to state the following theorem which we can prove using the last two estimate and the same proof of Theorem \ref{thm:Degen}.  
  
  \begin{thm} Let us consider the operator $Q_\alpha$ on $\mathbb{T}^d$ with $d\ge 2$.  If $\dim \sigma^* < \frac{\alpha^2}{4} \frac{s_d}{2} =  \frac{\alpha^2}{4} \frac{\pi^{d/2}}{\Gamma(d/2)}$, then $Q_\alpha \sigma =0$.
  Consequently, 
the operator $Q_\alpha$ is completely degenerate when $|\alpha|^2 >4d \frac{\Gamma(d/2)}{\pi^{d/2}}$.
  \end{thm}
  
  The about statement coincides with $d=1$ if we define $s_1=2$. 
  So,  the degeneracy condition $|\alpha|^2 >4d \frac{\Gamma(d/2)}{\pi^{d/2}}$ holds for all $d\ge 1$.
   
  \subsection{Decomposition of $H_\alpha$ and $Q_\alpha$, kernel and image of $\mathbb{E}Q_\alpha$ }
  
  The following result was known for $d=1$ (cf. \S \ref{subsect:decomp}). It will allows us to decompose the operator $Q_\alpha$ in a suitable way.
  Recall that the function $G$ is defined by (\ref{eq:G}).

  \begin{lem} \label{lem:decompd} There is a partition of $\mathbb{Z}^d\setminus \{0\}$ into two parts $\mathcal{P}_1$ and $\mathcal{P}_2$ such that
  the decomposition $ G = G_1 +G_2$ with 
  $$
     G_1(x) = \sum_{k\in \mathcal{P}_1} |k|^{-d}\exp(i k\cdot x), \qquad G_2(x) = \sum_{k\in \mathcal{P}_2} |k|^{-d}\exp(i k\cdot x)
  $$
  has the property
  $$
      G_1(x) =  \frac{1}{2} G(x) + O(1), \quad  G_1(x) =  \frac{1}{2} G(x) + O(1)  \ \ \ {\rm as }\ x \to 0.
  $$
  The partial sums, cubic or spheric, of both $G_1$ and $G_2$ are respectively bounded by $G_1(x) + C$
  or $G_2(x) +C$ for some constant $C$.   
  \end{lem}
  
  \begin{proof}
  We first look at the case  $d=2$. 
  Let 
  $$
      A = \begin{pmatrix} 
      1 & -1\\
      1&1
      \end{pmatrix} = \sqrt{2} \begin{pmatrix} 
      \cos \frac{\pi}{4} & -\sin \frac{\pi}{4}\\
      \sin \frac{\pi}{4}&\ \cos \frac{\pi}{4}
      \end{pmatrix}.
  $$ 
  The linear map defined by $A$ is a similitude such that $|Ak|= 2|k|$ and preserves the lattice $\mathbb{Z}^d$. 
  We choose $\mathcal{P}_1 = A(\mathbb{Z}^d\setminus\{0\})$ and  $\mathcal{P}_2 = \mathbb{Z}^d\setminus\{0\}\setminus \mathcal{P}_1$.
 Then
  $$
      G_1(x)  = \sum_{k\in \mathbb{Z}^d \setminus \{0\}} |Ak|^{-d}\exp(i Ak\cdot x)
      =\frac{1}{2}    \sum_{k\in \mathbb{Z}^d \setminus \{0\}} |k|^{-d}\exp(i k\cdot \, ^t\!A x).
  $$
  Now, by Theorem \ref{thm:G}, around $0$ we have
  $$
 G_1(x) =\frac{1}{2} G(A^tx) = \frac{1}{2} s_d \log \frac{1}{|^t\!Ax|} + O(1)  = \frac{1}{2} G(x) + O(1). 
 $$
 Consequently 
 $$
      G_2(x) = G(x) - G_1(x) =  \frac{1}{2} G(x) + O(1).
 $$
 When $d\ge 3$, we only need to replace the similitude $A$ by
 $$
      B = \begin{pmatrix} 
      A & \\
      & I_{d-2}
      \end{pmatrix} 
 $$
 where $I_{d-2}$ is the $(d-2) \times (d-2)$ unit matrix.  The boundedness of partial sums follows from Lemma  \ref{lem:partials}.
 \end{proof}
 
 Let $G_1^{(1)}=G_1$ and $G_2^{(1)}=G_2$.
 In the same way, we can continue to decompose $G_1^{(1)}(x)=G_1^{(2)}+ G_2^{(2)}$ with 
 $$
    G_1^{(2)}(x) = \frac{1}{2} G_1^{(1)}(^t\!A x)= \frac{1}{2^2} G( ^t\!A^2 x) =\frac{1}{4}   \sum_{k\in \mathbb{Z}^d \setminus \{0\}} |k|^{-d}\exp(i k\cdot\, ^t\!\!A^2 x).
 $$ 
 In our application later, we don't decompose $G_2^{(1)}$ or $G_2^{(2)}$ further. But we need to decompose the first component several times. 
 All these components $G_i^{(j)}$ have the same properties of $G$ as stated in Lemma \ref{lem:decompd}.

What we are really interested is $H_\alpha$, which is equal to $\frac{\alpha}{2}G(x)$.  The corresponding decomposition of $H_\alpha$ allows us to decompose 
the chaotic operator $Q_\alpha$ as follows
$$
   Q_{\alpha} = \widetilde{Q}_{\frac{\alpha}{\sqrt{2}^m}} Q_{\frac{\alpha}{\sqrt{2}^m}}\cdots Q_{\frac{\alpha}{\sqrt{2}^2}} Q_{\frac{\alpha}{\sqrt{2}}}
$$
where all the operators at the right are independent,  the correlation function of $Q_{\frac{\alpha}{\sqrt{2}^j}}$ has the same behavior as  $H_{\frac{\alpha}{\sqrt{2}^j}}$ and  
that of $\widetilde{Q}_{\frac{\alpha}{\sqrt{2}^m}}$ has the same behavior as  $H_{\frac{\alpha}{\sqrt{2}^m}}$.
\medskip

The kernel $\Phi_\alpha=e^{H_\alpha}$ is the Riesz kernel of order $\frac{\tau(d)\alpha^2}{2}$ (cf. (\ref{eq:Rieszd})).
If $\frac{\tau(d)\alpha^2}{2} <d$, the Lebesgue measure $\lambda$ is $Q_\alpha$-regular and the Peyri\`ere measure 
$\mathcal{Q}_\alpha$ associated to $Q_\alpha\lambda$ is well defined.  It can be proved that 
\begin{equation}\label{eq:LLN2d}
    \mathcal{Q}_\alpha\!-\!{\rm a.s.}  \qquad  \lim_{N\to \infty}   \frac{1}{\log N}  \sum_{n=1}^N\log P_n(t) =  \frac{\tau(d)\alpha^2}{4}.
         \end{equation}
         The proof is the same as that of  (\ref{eq:LLN2}).
         \medskip

 Now we can state the following description  for  the measure $Q_\alpha\sigma$
 when $\sigma$ is sufficiently regular.
 The proof is the same as that of Theorem \ref{thm:image}, but the above preparations are necessary.

 \begin{thm}\label{thm:imaged} Suppose that $D >\frac{\tau(d)\alpha^2}{4}$. If $\sigma \in \mathcal{R}_D^*$, then $Q_\alpha \sigma \in \mathcal{R}_{D-\frac{\tau(d)\alpha^2}{4}}^*$.
 \end{thm}
  
  Hence, following the same proof of Theorem \ref{thm:ker}, we can prove the following description of the kernel and image of the operator
  $\mathbb{E}Q_\lambda$.
  
  \begin{thm} \label{thm:kerd} Let $\tau(d) = \frac{\pi^{d/2}}{\Gamma(d/2)}$. We have 
  $$
  \mathcal{S}_{\frac{\tau(d) \alpha^2}{4} -} \subset {\rm Ker} \mathbb{E} Q_\alpha \subset \mathcal{S}_{\frac{\tau(d) \alpha^2}{4} +};
  $$
   $$
  \mathcal{R}_{\frac{\tau(d) \alpha^2}{4} +} \subset {\rm Im} \ \mathbb{E} Q_\alpha \subset \mathcal{R}_{\frac{\tau(d) \alpha^2}{4} -}.
  $$
  \end{thm}
 
 Let us finish by stating the counterpart in higher dimension of Theorem \ref{thm:dim}
 \begin{thm} \label{thm:dimd}Suppose that  $\sigma$ is a unidimensional measure on $\mathbb{T}^d$ and $\dim \sigma >\frac{\tau(d)\alpha^2}{4}$. We have $$\dim Q_\alpha \sigma = \dim \sigma - \frac{\tau(d)\alpha^2}{4}.$$
\end{thm}

 \section{Random trigonometric series} 
 
 \subsection{Study of $\sum_{n=1}^\infty \rho_n \cos  (nt +\omega_n)$}
 Now come back to our series
 \begin{equation}\label{eq:RS*}
 \sum_{n=1}^\infty X_n(t) = \sum_{n=1}^\infty \rho_n \cos  (nt +\omega_n),
\end{equation} 
which was proposed in \S 1 Introduction. Suppose that $\sum \rho_n^2 =\infty$.
This  series diverges a.s. almost everywhere with respect to the Lebesgue measure. Indeed, for any $t$ fixed,
the series diverges almost surely (Theorem 4, page 31 in \cite{Kahane1985} applies). Then Fubini theorem allows us to conclude.
The same argument shows that for any given measure $\sigma$ on $\mathbb{T}$, the series 
diverges almost surely $\sigma$-almost everywhere.  
 When $\rho_n=\frac{1}{\sqrt{n}}$, almost surely the series diverges everywhere (cf. \cite{Kahane1985}, p. 108).
In the following we investigate how the series diverges. As we shall see, we can obtain  speeds of divergence 
for points in the supports of our chaotic measures $Q_\alpha\lambda$. Similar results also hold for $Q_\alpha\sigma$ with $\sigma \not=\lambda$, but we don't 
state them.


Let $\{\alpha_n\}$ be a sequence of real numbers 
which defines a chaotic operator $Q$.
Suppose that the Lebesgue measure $\lambda$ is $Q$-regular. For example,   $\lambda$ is $Q_\alpha$-regular
if $|\alpha|<2$ (cf. Theorem \ref{thm:image}).
 Consider  the 
Peyri\`ere measure $\mathcal{Q}$ associated to $Q\lambda$. By Proposition \ref{prop:moments}, we have
     \begin{eqnarray*}
     \mathbb{E}_{\mathcal{Q}} X_n(t) 
    & = & 
    \frac{1}{2}\cdot \rho_n\alpha_n+ O(\rho_n |\alpha_n|^3)\\
     \mathbb{E}_{\mathcal{Q}} X_n(t)^2 
    & = & 
      \frac{1}{2} \cdot \rho_n^2+ O(\rho_n^2 |\alpha_n|^2)
   \end{eqnarray*}
   Remark that the variance of $X_n$ is equivalent to $\frac{1}{2}\rho_n^2$ if $\alpha_n \to 0$. Notice that we have coefficients $\rho_n$ in the series (\ref{eq:RS*}), but 
   our chaotic measures are defined by the exponentiation of $\alpha_n \cos (nt+\omega_n)$, where $\alpha_n$'s are not necessarily $\rho_n$'s.
   
   \begin{thm}\label{thm:ILL0} Suppose $\sum \rho_n^2 =\infty$ and $\lambda$ is  $Q$-regular where $Q$ is the 
   chaotic operator defined by $\{\alpha_n\}$ with $\alpha_n \to 0$. Then almost surely $Q \lambda$-almost everywhere the series (\ref{eq:RS*}) diverges  and 
   $$
       \sum_{n=1}^N \rho_n \cos (n x +\omega_n) = \frac{1}{2} \sum_{n=1}^N \rho_n\alpha_n+ O\left(\varphi\Big(\sum_{n=1}^N \rho_n^2\Big)\right),
   $$
  where $\varphi : [0, \infty)\to [0, \infty)$ is a function satisfying the condition
  $$
       \sum_{n=1}^\infty \frac{\rho_n^2}{\varphi\Big(\sum_{n=1}^N \rho_n^2\Big)^2}<\infty.
  $$ 
   \end{thm}
   \begin{proof} The divergence is already discussed  at the beginning of this section, because $\{nt +\omega_n)\}$ are 
   $\mathcal{Q}$-independent. 
   The proof for the estimation is the same as that in Proposition \ref{prop:LogP}. Instead we consider the series
   $$
         \sum_{n=1}^\infty \frac{X_n(t) - \mathbb{E}_{\mathcal{Q}_\alpha} X_n(t)}{\varphi\Big(\sum_{n=1}^N \rho_n^2\Big)},
   $$
   the partial sums of which form a $L^2$-bounded martingale.
   This boundedness is checked by the fact that
   the variance of $X_n(t)$ is of size $\rho_n^2$.
   \end{proof}
   
   We can improve the above result by a  law of iterated logarithm as a special case of the   law of iterated logarithm obtained  by Wittmann \cite{Wittmann1985}. 
   
      \begin{thm}\label{thm:ILL} Suppose $\sum \rho_n^2 =\infty$  and $\lambda$ is  $Q$-regular where $Q$ is the 
   chaotic operator defined by $\{\alpha_n\}$ with $\alpha_n \to 0$. Then almost surely $Q \lambda$-almost everywhere  we have
   $$
       \limsup_{n\to \infty} \frac{\sum_{n=1}^N \rho_n \big( \cos (n x +\omega_n) -  \frac{1}{2} \alpha_n\big)}
       {{\sqrt{\sum_{n=1}^N \rho^2_n \log \log \sum_{n=1}^N \rho^2_n} } } =1
   $$
   if the following condition is satisfied for some $2<p\le 3$: 
   \begin{equation}\label{eq:condILL}
       \sum_{n=1}^\infty  \frac{\rho_n^p}{ \Big(\sqrt{\sum_{n=1}^N \rho^2_n  \log \log \sum_{n=1}^N \rho^2_n} \Big)^{p} } <\infty.
   \end{equation}
   
   \end{thm}
   \begin{proof} It is a special case of the following law of iterated logarithm due to Wittmann (cf. \cite{Wittmann1985} Theorem 1.2). Let $\{Z_n\}$ be a sequence of independent real random variables
   such that $\mathbb{E}Z_n=0$ and $\mathbb{E}Z_n^2<\infty$. Let 
   $$
      S_n = \sum_{k=1}^n Z_k, \qquad s_n = \sqrt{\mathbb{E} S_n^2}, \qquad t_n = \sqrt{2 \log \log s_n^2}.
   $$
   Suppose that    
   \begin{equation}\label{eq:wittman}
       \lim s_n =\infty, \quad s_{n+1}=O(s_n), \quad  \sum_{n=1}^\infty \frac{\mathbb{E} |Z_n|^p}{(s_nt_n)^p} <\infty \ \ (\exists 2<p\le 3).
   \end{equation}
   Then almost surely we have 
   $\limsup_n \frac{S_n}{s_nt_n}=1$.
   We can apply this result to $Z_n = X_n(t) - \mathbb{E}_{\mathcal{Q}} X_n(t)$ where $X_n(t) = \rho_n \cos (nt +\omega_n)$ to obtain the announced result. 
   Indeed, 
   we have 
   $$
       s_n^2 \sim \frac{1}{2} \sum_{k=1}^n \rho_n^2.
   $$
   So, since $|Z_n| \le 2 \rho_n$, the convergence of the series in (\ref{eq:wittman}) is ensured by the condition (\ref{eq:condILL}). 
   The assumption $\sum \rho_n^2=\infty$ is nothing but $\lim s_n =\infty$. The condition 
    $s_{n+1}=O(s_n)$ in (\ref{eq:wittman})  is also satisfied because $\rho_n \to 0$.
   \end{proof}

   In Theorem  \ref{thm:ILL0} and Theorem \ref{thm:ILL}, we can consider the operator $Q$ as well as  $Q'$ associated to $\{\eta_n\alpha_n\}$
   where $\{\eta_n\}$ is a sequence of $+1$ or $-1$ and then consider  the measures $Q \lambda$ and   $Q'\lambda$.  
   Both operators $Q$ and $Q'$ have the same correlation function. 
      Both  measures $Q'\lambda$ and $Q\lambda$
   share same properties,  like $$
   \dim Q'_{\alpha}\lambda = \dim Q_{\alpha}\lambda= 1-\frac{\alpha^2}{4}\ \  (|\alpha|<2)
   $$
   in the case $\alpha_n = \frac{\alpha}{\sqrt{n}}$.
   But they are usually mutually singular (cf. Proposition \ref{prop:mutualsing}):
   $$ 
      a.s \ \ \  Q' \lambda \perp  Q \lambda \Longleftrightarrow  \sum_{n=1}^\infty |1-\eta_n|^2 |\alpha_n|^2=\infty.
   $$
    By choosing different sequences $\{\alpha_n\}$ and $\{\eta_n\}$, points of different properties  the series (\ref{eq:RS*})  can be obtained, even when $\{\alpha_n\}$ is fixed. 
    Through all these chaotic measures, we see that the partial sums of the series (\ref{eq:RS*})
    are very multifractal.

 \subsection{Special series $\sum_{n=1}^\infty \frac{\cos  (nt +\omega_n)}{n^r}$ ($0<r \le \frac{1}{2}$) }\, \ \ \, 
   Let us look at the special series 
 \begin{equation}\label{eq:series1/2a}
   \sum_{n=1}^\infty \frac{\cos  (nt +\omega_n)}{n^r} \qquad (0<r\le \frac{1}{2}).
 \end{equation}
 
 First assume $r=\frac{1}{2}$. The variables $nt+\omega_n$ are i.i.d. with respect to $\lambda \otimes P$ and we have
 $$
       \mathbb{E}_{\lambda\otimes P} \cos  (nt+\omega_n)=0,
       \quad   \mathbb{E}_{\lambda\otimes P} \cos^2  (nt+\omega_n)=\int \cos^2  x dx = \frac{1}{2}.
 $$
 Then, by the classical law of iterated logarithm, a.s. for $\lambda$-almost all $t$ we have
 $$
     \limsup_{N\to \infty} \frac{\sum_{n=1}^N \cos  (nt +\omega_n)}{\sqrt{N\log\log  N}} = 1.
 $$
 By the law of iterated logarithm in
  \cite{Wittmann1985},   we have a.s. Lebesgue-almost everywhere
 $$
    \limsup_{N\to \infty} \frac{\sum_{n=1}^N \frac{\cos  (nt +\omega_n)}{\sqrt{n}}}{\sqrt{ \log N \log\log \log N}} =1.
 $$
 More generally, if $|\alpha|<2$, by Theorem  \ref{thm:ILL},  we have a.s. $Q_\alpha\lambda$-almost everywhere
  \begin{equation}\label{eq:ILL12}
    \limsup_{N\to \infty} \frac{\sum_{n=1}^N \frac{1}{\sqrt{n}}  \cos  (nt +\omega_n)-\frac{\alpha}{2}\log N}{\sqrt{ \log N \log\log \log N}} =1.
 \end{equation}
 In particular,  we have a.s. $Q_\alpha\lambda$-almost everywhere
  \begin{equation}\label{eq:LLN12}
    \lim_{N\to \infty} \frac{1}{\log N} \sum_{n=1}^N \frac{1}{\sqrt{n}}  \cos  (nt +\omega_n)=\frac{\alpha}{2}.
 \end{equation}

 Now assume $0\le r<\frac{1}{2}$. Since $\sum_{k=1}^n k^{-1/2-r}= \frac{n^{1/2-r}}{1/2-r}+O(1)$ and $\sum_{k=1}^n k^{-2r}= \frac{n^{1-2r}}{1-2r}+O(1)$, 
  if $|\alpha|<2$, by Theorem  \ref{thm:ILL},  we have a.s. $Q_\alpha\lambda$-almost everywhere
    \begin{equation*}
    \limsup_{N\to \infty} \frac{\sum_{n=1}^N \frac{1}{n^r}  \cos  (nt +\omega_n)-\frac{\alpha}{1-2r} N^{1/2-r}}{\sqrt{ N^{1-2r} \log \log N}} =\sqrt{1-2r}.
  \end{equation*}
  Since  $N^{1/2 -r}=o(\sqrt{ N^{1-2r} \log \log N})$,   we have a.s. $Q_\alpha\lambda$-almost everywhere
   \begin{equation}\label{eq:ILLr}
    \limsup_{N\to \infty} \frac{\sum_{n=1}^N \frac{1}{n^r}  \cos  (nt +\omega_n)}{\sqrt{ N^{1-2r} \log \log N}} =\sqrt{1-2r}.
  \end{equation}
  \bigskip
  
  It is well known that the series  $\sum_{n=1}^\infty \frac{\sin  nt}{\sqrt{n}}$ converges everywhere, but not uniform, to a Lebesgue integrable function, 
  and the convergence is uniform on any interval $[\delta, 1-\delta]$ ($0<\delta<1/2$)
  (cf. \cite{Zygmund}, vol. I, Chapter V and \cite{Bary1964} p.87).   But the randomized series 
   $\sum_{n=1}^\infty \frac{\sin  (nt +\omega_n)}{\sqrt{n}}$ loses  all these properties almost surely. 

\subsection{Large deviation of $\sum_{n=1}^N \frac{\cos (nt+\omega_n)}{\sqrt{n}}$}

Suppose $|\alpha|<2$. By (\ref{eq:ILL12}) which is a consequence of Theorem \ref{thm:ILL0}, we have 
$$
   \mathcal{Q}_\alpha\!-\!a.s.  \quad 
   \lim_{N\to \infty}\frac{1}{\log N} \sum_{n=1}^N \frac{\cos (nt+\omega_n)}{\sqrt{n}} =\frac{\alpha}{2}.
$$
We have the following stronger result, a large deviation result.

\begin{thm} \label{thm:LD}Suppose $|\alpha|<2$. For any $\eta >0$,  we have
$$
    \lim_{N\to \infty} \frac{1}{\log N} \log \mathcal{Q}_\alpha\left\{ (t, \omega):  \frac{1}{\log N}\sum_{n=1}^N \frac{\cos (nt+\omega_n)}{\sqrt{n}}  \not\in \frac{\alpha}{2}+ [-\eta, \eta]\right\}
    =-\eta^2.
$$
\end{thm}

\begin{proof}
Let $W_n = \sum_{k=1}^n \frac{\cos (kt+\omega_k)}{\sqrt{k}}$, the random variable in question and   $a_n=\log n$, the chosen normalizer.   We shall prove that  the following limit exits
$$
   c(\beta):= \lim_{n\to \infty} \frac{1}{a_n} \log \mathbb{E}_{\mathcal{Q}_\alpha} e^{\beta W_n}= \frac{(\beta+\alpha)^2 -\alpha^2}{4}, \quad \forall \beta \in \mathbb{R}
$$
which is called  the { free energy function} of $(W_n)$  with weight $(a_n)$, with respect to the probability measure $\mathcal{Q}_\alpha$. 
By the large deviation theorem (\cite{Ellis1985}, p. 230), for any non empty interval $K\subset \mathbb{R}$ we have
$$
    \lim_{n\to \infty} \frac{1}{a_n} \log \mathcal{Q}_\alpha\big\{ a_n^{-1} W_n \not\in K\big\}
    = - \inf_{\gamma \in K} c^*(\gamma)
$$ 
where  $c^*(\gamma)=\sup_\beta (\gamma \beta -c(\beta))$ is the Legendre transform of $c(\cdot)$. 
Then the announced result will follow. 

Indeed, the limit $c(\beta)= \frac{\beta^2 +2\alpha \beta}{4}$ is easy to obtain, because 
$$
    \mathbb{E}_{\mathcal{Q}_\alpha} e^{\beta W_n}= \prod_{k=1}^n \frac{I_0\big(\frac{\alpha +\beta}{\sqrt{k}}\big)}{I_0\big(\frac{\alpha}{\sqrt{k}}\big)}
    \asymp \exp \left(\frac{\beta^2 +2\alpha \beta}{4} \log n \right).
$$
It is also direct to get the Legendre transform $c^*(\gamma) = \big(\gamma-\frac{\alpha}{2}\big)^2$.  The convex function 
$c^*(\cdot) $ attains its minimal value at $\gamma = \frac{\alpha}{2}$. 
For $K = \frac{\alpha}{2}+ [-\eta, \eta]$,  
we have clearly 
$$\inf_{\gamma\in K} c^*(\gamma) = c^*(\alpha/2 +\eta)=  \eta^2.$$ 
\end{proof}

The same method applies to other cases. Assume $\rho_n \to 0$ and $\frac{1}{\sqrt{n}} = o(\rho_n)$. That is the case for $\rho_n =\frac{1}{n^r}$ with $0<r<\frac{1}{2}$. If we consider
$W_n = \sum_{k=1}^n \rho_n  \cos (kt+\omega_k)$, we have
$$
    \mathbb{E}_{\mathcal{Q}_\alpha} e^{\beta W_n}= \prod_{k=1}^n \frac{I_0\big(\beta \rho_k + \frac{\alpha}{\sqrt{k}}\big)}{I_0\big(\frac{\alpha}{\sqrt{k}}\big)}
    \asymp \exp \left(\frac{\beta^2}{4} \sum_{k=1}^n \rho_k^2 \right).
$$
Take $a_n= \sum_{k=1}^n \rho_k^2$ as normalizer. Then we get that the  free energy function of $(W_n)$ is equal to
$c(\beta)=\frac{\beta^2}{4}$. Its Legendre transform $c^*(\gamma)= \gamma^2$. Therefore
$$
    \lim_{N\to \infty} \frac{1}{\sum_{n=1}^N \rho_n^2} \log \mathcal{Q}_\alpha \left\{ (t, \omega):  \frac{\sum_{n=1}^N \rho_n \cos (nt+\omega_n)}{\sum_{n=1}^N \rho_n^2}\not\in [-\eta, \eta] \right\}
    =-\eta^2.
$$

\end{document}